\documentclass[11pt]{article}
\usepackage{geometry}
\usepackage{amsmath,amssymb,amsthm}
\usepackage{charter}
\usepackage{cite}
\usepackage{xcolor}
\usepackage[nottoc]{tocbibind} 
\usepackage[pagebackref=true]{hyperref}
\usepackage{comment}
\renewcommand*\backref[1]{\ifx#1\relax \else \mbox{\textcolor{gray}{Page #1}} \fi}
\hypersetup{
colorlinks,%
citecolor=blue,%
filecolor=black,%
urlcolor=black,%
linkcolor=blue
}

\newcommand\blfootnote[1]{%
\begingroup
\renewcommand\thefootnote{}\footnote{#1}%
\addtocounter{footnote}{-1}%
\endgroup
}

\newtheorem{teo}{Theorem}
\newtheorem{prop}[teo]{Proposition}
\newtheorem{lema}[teo]{Lemma}
\newtheorem{cor}[teo]{Corollary}
\theoremstyle{definition} 

\DeclareMathOperator{\V}{SL}
\newcommand{\N}{_n(\mathbb{R})}

\newcommand{\dif}{\mathop{}\!\mathrm{d}}

\newcommand{\RR}{(\mathbb{R}^n)}
\newcommand{\NN}{(\mathbb{R})}
\DeclareMathOperator{\conv}{conv}
\DeclareMathOperator{\Conv}{Conv}

\DeclareMathOperator{\I}{I}

\DeclareMathOperator{\PA}{Conv_{p.a}}
\DeclareMathOperator{\LQ}{_{l.q}}

\DeclareMathOperator{\dom}{dom}
\DeclareMathOperator{\interior}{int}

\DeclareMathOperator{\MA}{MA}
\DeclareMathOperator{\MAh}{MA,h}
\DeclareMathOperator{\MAp}{MA,p}

\DeclareMathOperator{\D}{D^2}

\DeclareMathOperator{\diag}{diag}
\DeclareMathOperator{\epi}{epi}
\DeclareMathOperator{\LC}{Conv_{ld}}
\DeclareMathOperator{\LP}{Conv_{lp}}

\DeclareMathOperator{\bd}{bd}
\DeclareMathOperator{\as}{\Omega}
\DeclareMathOperator{\fas}{\Omega_{func}}

\newcommand{\finite}{(\mathbb{R}^n; \mathbb{R})}

\DeclareMathOperator{\dist}{dist}
\DeclareMathOperator{\diam}{diam}

\DeclareMathOperator{\supp}{supp}
\DeclareMathOperator{\conc}{Conc([0,\infty))}

\newcommand{\seq}{\xrightarrow{\tau}}

\makeatletter
\renewcommand{\@fnsymbol}[1]{}
\makeatother

\begin{document}
\title{A Characterization of Functional Affine Surface Areas} \author{Fernanda M. Baêta\thanks{Departamento de Matemática, ICEx, Universidade Federal de Minas Gerais, 31270-901 Belo Horizonte, Brazil. Email: fernandah@ufmg.br}}

\date{}

\maketitle
\pretolerance10000
\begin{abstract}
A characterization of valuations on the space of convex Lipschitz  functions whose domain is a polytope in $\mathbb{R}^n$ is obtained. It is shown that every upper semicontinuous, equi-affine and dually epi-translation invariant valuation can be written as a linear combination of a constant term, the volume of the domain, and a functional affine surface area. In addition, dual statements for finite-valued convex functions are established.

\end{abstract}

\section{Introduction and Statement of the Main Result}

Let $K\subset \mathbb{R}^{n}$ be a convex body (i.e., a non-empty compact convex set) and denote by $\bd K$ its topological boundary. By a theorem of Aleksandrov (see \cite{aleksandrov1939second}), for almost every point $x\in \bd K$ with respect to the $(n-1)$-dimensional Hausdorff measure $\mathcal{H}^{n-1}$, there exists a paraboloid that osculates $\bd K$ at $x$. By Schütt and Werner \cite{schuett1990convex}, the \textit{affine surface area} of a general (not necessarily smooth) convex body is defined by
\begin{align*}
\as(K)= \int_{\bd K}\kappa(x)^\frac{1}{n+1}\dif \mathcal{H}^{n-1}(x),    
\end{align*}
where the (generalized) Gaussian curvature $\kappa(x)$ of $\bd K$ at such a point  $x$ is taken as the Gaussian curvature of the corresponding osculating paraboloid. For smooth convex surfaces, this definition coincides with the classical one (see \cite{blaschke1923differentialgeometrie}).  

The affine surface area $\as$ is defined for every $K \in \mathcal{K}^n$, where $\mathcal{K}^n$ denotes the space of convex bodies equipped with the topology induced by the Hausdorff metric. Lutwak \cite{lutwak1991extended} showed that $\as$ is upper semicontinuous, that is, if $K_m \to K$ in the Hausdorff metric, then 
$$ \as(K)\ge \limsup_{m \to \infty} \as(K_m).$$ 
Moreover,
$\as$ is \textit{equi-affine invariant}, i.e., invariant under volume-preserving affine maps, and 
a valuation (see \cite{schutt1993affine}). More precisely, a functional $\Phi: \mathcal{K}^n\rightarrow \mathbb{R}$ is called a \textit{valuation} if 
\begin{align*}
    \Phi(K\cup L)+ \Phi(K\cap L)= \Phi(K)+\Phi(L)
\end{align*}
whenever $K,L$ and  $K\cup L$ belong to $\mathcal{K}^n$.

One of the most famous results 
on valuations is the Hadwiger characterization theorem, which classifies all continuous and rigid motion invariant valuations on the space $\mathcal{K}^n$ (see \cite{hadwiger1957vorlesungen}). Specifically, it states that any such valuation can be expressed as a linear combination of the intrinsic volumes $V_0, V_1, \dots, V_n$,  which include the $n$-dimensional volume $V_n$, the Euler characteristic $V_0$, and the surface area $2V_{n-1}$.
A complete characterization of affine surface area is due to Ludwig and Reitzner, representing an equi-affine 
analog of Hadwiger's theorem.

\begin{teo}[\cite{ludwig1999affine,ludwig1999characterization}]\label{cg}
A functional  $\Phi: \mathcal{K}^n\rightarrow \mathbb{R}$ is an upper semicontinuous and equi-affine invariant valuation if and only if there exist constants $c_0,c_1$ and $c_2\geq 0$ such that
\begin{align*}
    \Phi(K)= c_0V_0(K)+c_1\,V_n(K)+c_2\Omega(K)
\end{align*}
for every $K\in \mathcal{K}^n$.
\end{teo}

Ludwig and Reitzner \cite{ludwig2010classification} 
further posed the problem of classifying all valuations on the space $\mathcal{K}_0^n$  of convex bodies containing the origin in their interiors, which are invariant under volume-preserving linear transformations. They showed that any valuation $\Phi: \mathcal{K}_0^n\rightarrow \mathbb{R}$ that is  upper semicontinuous, invariant under volume preserving linear transformations, and  vanishes on polytopes can be represented as
\begin{align*}
    \Phi(K)= \int_{\bd K} \zeta(\kappa_0(x))\dif \mu_K(x),
\end{align*}
where $\dif \mu_K(x)=\langle x,n(x)\rangle \dif \mathcal{H}^{n-1}(x)$ is the cone measure on $\bd K, n(x)$ is the exterior  unit normal vector to $K$ at $x\in\bd K$,  defined for $\mu_K$-a.e.\ $x$, and 
$$\kappa_0(x)=\frac{\kappa(x)}{\langle x,n(x)\rangle^{n+1}}.$$
Here, $\zeta$ belongs to the set
\begin{align}\label{concave}
\conc=\left\{\zeta: [0,\infty)\rightarrow [0,\infty)\mid \ \zeta \ \mbox{is concave}, \ \lim_{t\rightarrow 0}\zeta(t)=0,  \mbox{ and } \lim_{t\rightarrow \infty}\zeta(t)/t =0\right\}   
\end{align}
(see \cite{ludwig2010classification, ludwig2001semicontinuity}), and   $\langle x,y \rangle$ denotes  the standard inner product of $x$ and $y$, either in $\mathbb{R}^n$ or  $\mathbb{R}^{n+1}$ depending on the context. A complete answer to Ludwig and Reitzner’s problem in the upper semicontinuous case was later provided by Haberl and Parapatits (see \cite{haberl2014centro}).  Further characterizations of valuations with interesting invariance properties were obtained, for example, in \cite{alesker,haberl, haberl2, hadwiger1957vorlesungen,klain,L, mullen,mullen2}. For results on the affine surface area of convex bodies, we refer to \cite{ludwig1999characterization, lutwak1991extended, schuett1990convex}; see also the recent survey \cite{schuett2023affine}.

One important direction of generalization is the extension of valuations from convex bodies to convex functions, which has attracted significant attention in recent years; see, for instance, \cite{colesanti2019convex, colesanti2017minkowski, colesanti2020homogeneous, colesanti2020hessian}. 
Further developments and characterizations can be found in \cite{M1, colesanti2023hadwiger, colesanti2024hadwigerI, knoerr2021support, knoerr2024valuations}, and more recently in \cite{mouamine2025klain, mouamine2025vectorial, mouamine2025additive}.
Denote by  $\Conv(\mathbb{R}^n)$ the space of  lower semicontinuous, convex, and proper functions (i.e., not identically $\infty$) defined on $\mathbb{R}^n$, and let
\begin{align*}
\Conv(\mathbb{R}^n;\mathbb{R})= \{v:\mathbb{R}^n\to \mathbb{R}\mid \ v \text{ is convex}\}    
\end{align*}
be the space of finite-valued convex functions on $\mathbb{R}^n$. A functional $Z: \Conv(\mathbb{R}^n)\rightarrow \mathbb{R}$ is called a \textit{valuation} if
$$Z(u\wedge v)+ Z(u\vee v)= Z(u) + Z(v)$$
for every $u,v\in \Conv(\mathbb{R}^n)$ such that also $u\wedge v, u\vee v\in \Conv(\mathbb{R}^n)$, where $u\wedge v$ and $u\vee v$ denote the
pointwise  minimum and maximum of $u,v\in \Conv(\mathbb{R}^n)$, respectively.

Let $u: U\subset \mathbb{R}^n \to \mathbb{R}$ be a convex function, where $U$ is an open bounded convex set, and let
\(
\epi(u)\subset \mathbb{R}^{n+1}
\)
denote its epigraph. Then the affine surface area of $\epi(u)$ is given by
\begin{align*}
    \fas(u)=\int_U (\det \D u(x))^{\frac{1}{n+2}} \dif x,
\end{align*}
where $\D u(x)$ denotes the Hessian matrix of $u$ at $x$ and $\D u(x)=0$ if $u$ is not twice differentiable at $x$ (see \cite{trudinger2000bernstein}). According to Aleksandrov's theorem \cite{aleksandrov1939second}, a convex function
$u:U\subset\mathbb{R}^n\to\mathbb{R}$ is twice differentiable almost everywhere. In particular, the Hessian $\D u(x)$ exists for almost every $x\in U$. Therefore, the integral above is well-defined (but possibly infinite).
For further definitions of functional affine surface area, see, for example, \cite{artstein2012functional, caglar2016functional, li2019floating, schutt2025weighted}.

The set of points where the function $u\in \Conv(\mathbb{R}^n)$ takes finite values is called the (effective) \emph{domain} of $u$, denoted by 
$$
\dom u = \{ x \in \mathbb{R}^n \mid \ u(x) < \infty \},
$$  
and for a set $A \subseteq \mathbb{R}^n$, we denote its interior by $\interior(A)$.  
 We  consider the subspace of $\Conv(\mathbb{R}^n)$,
\begin{align*}
 \LC\RR=
\{u\in \Conv(\mathbb{R}^n)\mid \ \dom u \mbox{ is compact},\, u \mbox{ is Lipschitz on  $\interior(\dom u)$}\},
\end{align*}
that is, the set of functions $u \in \Conv(\mathbb{R}^n)$ with compact domain that are Lipschitz on the interior of their domain. In particular,  if the interior of $\dom u$ is empty, then $u$ is not required to be Lipschitz. 

By \cite[Theorem~1]{baetaludwigsemicontinuity}, the functional $\fas$ is finite on $\LC(\mathbb{R}^n)$. Moreover, 
the functional $\fas: \LC(\mathbb{R}^n) \to \mathbb{R}$ is a valuation, and it satisfies the following properties.

\begin{itemize}
    \item[(i)]   $\V\N$ \textit{invariance}, i.e., $\fas(u \circ \varphi) = \fas(u)$ for every $u \in \LC\RR$ and $\varphi \in \V\N$;
 
    \item[(ii)]   \textit{translation invariance}, i.e., $\fas(u \circ \tau) = \fas(u)$
    for every $u \in \LC\RR$ and translation $\tau$ on $\mathbb{R}^n$;

    \item[(iii)]  \textit{vertical translation invariance}, i.e., $\fas(u+c) = \fas(u)$ for every $u \in \LC\RR$ and every constant $c \in \mathbb{R}$;
    
    \item[(iv)]  \textit{dual translation invariance}, i.e., $\fas(u+l) = \fas(u)$ for every every $u \in \LC\RR$ and every  linear function $l$ on $\mathbb{R}^n$. 
\end{itemize}
A valuation $Z: \Conv(\mathbb{R}^n)\to \mathbb{R}$ that is $\V\N$-invariant and translation invariant is called an \textit{equi-affine invariant valuation}.  
If $Z$ is both vertically and dually translation invariant, we say that $Z$ is \textit{dually epi-translation invariant}.  
The same definitions apply when $Z$ is defined on a subspace of $\Conv(\mathbb{R}^n)$.

We equip the set $\LC(\mathbb{R}^n)$ with the following notion of convergence: We say that $u_k\in \LC(\mathbb{R}^n)$, with $k\in\mathbb{N}$, is $\tau$-convergent to $u\in \LC(\mathbb{R}^n)$, written $u_k \seq u$, if for every $x\in \mathbb{R}^n$ the following hold:
\begin{enumerate}
\item[(i)] for every sequence $x_k \to x$, we have 
$$u(x)\leq \liminf_{k\rightarrow \infty} u_k(x_k);$$
\item[(ii)] there exists a sequence $x_k \to x$ such that 
$$u(x)= \lim_{k\rightarrow \infty}u_k(x_k);$$
\item[(iii)] the Lipschitz constants of $u_k$ on the interior of $\dom u_k$ are uniformly bounded by some constant $M>0$ independent of $k$.
\end{enumerate}

\noindent This notion of convergence was also considered in \cite{baeta1d,baetaludwigsemicontinuity}. A related formulation of condition (iii) was first considered in \cite{colesanti2020invariant, colesanti2021continuous}, in the context of real-valued Lipschitz continuous maps on the unit sphere $S^{n-1}$. The convergence given by conditions (i) and (ii) is called \textit{epi-convergence} (for more details see  \cite{rockafellar2009variational}).

We say that  $Z:\LC\RR\rightarrow \mathbb{R}$ is  $\tau$-\textit{upper semicontinuous}  if, for every sequence $u_k$ in $\LC\RR$, $k\in\mathbb{N}$,  that is  $\tau$-convergent   to $u\in \LC\RR$,
$$Z(u)\geq \limsup_{k\rightarrow \infty}Z(u_k).$$

Note that the function $t \mapsto t^{\frac{1}{n+2}}$, for $t \ge 0$, belongs to the class $\conc$ defined in \eqref{concave}, and is precisely the function appearing in the definition of $\fas$. Moreover, the restriction of $\fas$ to $\LC(\mathbb{R}^n)$ is $\tau$-upper semicontinuous, as proved in \cite[Theorem~1]{baetaludwigsemicontinuity}. We consider the following subset of $\LC(\mathbb{R}^n)$,
\begin{align*}
   \LP(\mathbb{R}^n) = \{\, u \in \LC(\mathbb{R}^n) \mid \ \dom u \text{ is a polytope} \,\}.
\end{align*}

Our goal is to prove a functional version of Theorem~\ref{cg}, thereby obtaining a characterization of the more general functional affine surface area
\begin{align}\label{int}
Z(u) = \int_{\dom u} \zeta(\det \D u(x))\dif x
\end{align}
for $\zeta \in \conc$ and $u \in \LP\RR$. This functional version is naturally related to the geometric version in $\mathbb{R}^{n+1}$, since the epi-graph of a convex function $u \in \LP\RR$ is a convex subset of $\mathbb{R}^{n+1}$. To the best of our knowledge, Theorem~\ref{equiv} provides the first characterization of upper semicontinuous valuations on the space of convex functions in arbitrary dimension $n$. The one-dimensional case was established previously (see~\cite{baeta1d}), and our result extends this to the general setting.

\begin{teo}\label{equiv}
A functional $Z: \LP\RR\rightarrow\mathbb{R}$ is a $\tau$-upper semicontinuous, equi-affine  and dually epi-translation invariant valuation if and only if there are  constants $c_0,c_1\in\mathbb{R}$ and a function $\zeta\in\conc$ such that
\begin{align}\label{maintheo}
Z(u)= c_0+c_1\,V_n(\dom u) +\int_{\dom u}\zeta(\det \D u(x))\dif x
\end{align}
for every $u\in \LP\RR$.
\end{teo}  
By \cite[Theorem 1]{baetaludwigsemicontinuity} and \cite[Lemma 12]{baeta1d}, the functional \eqref{int} defines an equi-affine and dually epi-translation invariant valuation that is finite and $\tau$-upper semicontinuous for all $u\in\LP(\mathbb{R}^n)$. Furthermore, by \cite[Lemma~6 and Lemma~9]{baeta1d}, the volume of the domain is a $\tau$-continuous, equi-affine and dually epi-translation invariant valuation. Therefore, the functional on the right-hand side of \eqref{maintheo} is a $\tau$-upper semicontinuous, equi-affine and dually epi-translation invariant valuation. We expect that additional $\tau$-upper semicontinuous, equi-affine and dually epi-translation invariant valuations exist on $\LC(\mathbb{R}^n)$ which vanish on $\LP(\mathbb{R}^n)$. One such candidate is given by the affine surface area of the domain. This is consistent with the classical fact that affine surface area vanishes on polytopes, and is supported by the fact that if $u_k\in \LC(\mathbb{R}^n), k\in\mathbb{N}, $ is $\tau$-convergent to $u\in \LC(\mathbb{R}^n)$ and each $\dom u_k$ has non-empty interior, then $\dom u_k$ converges to $\dom u$ in the Hausdorff metric (see \cite[Lemma~9]{baeta1d}).

The uniform boundedness of the Lipschitz constants in the definition of $\tau$-convergence is essential. Indeed, \cite[Example~1]{baetaludwigsemicontinuity} shows that, without this assumption, the general functional affine surface area \eqref{int} is, in general, no longer $\tau$-upper semicontinuous. It is worth noting, moreover, that there exist sequences of functions with compact domains that epi-converge to a function with a compact domain, yet the volumes of the domains do not converge to the volume of the domain of the limit function (see also \cite[Example~1]{baetaludwigsemicontinuity}).

Note  that the integral \eqref{int}  is not $\tau$-continuous, as it vanishes on piecewise affine functions. As a consequence of Theorem \ref{equiv}, we obtain the following result.

\begin{cor}\label{cor1}
A functional  $Z: \LC\RR\rightarrow\mathbb{R}$ is a $\tau$-continuous, equi-affine  and dually epi-translation invariant valuation if and only if there are constants $c_0,c_1\in\mathbb{R}$ such that
$$Z(u)=c_0+c_1\,V_n(\dom u)$$ 
for every $u\in \LC\RR$.
\end{cor}

Note that this result is established in $\LC\RR$, since every $u \in \LC\RR$ can be approximated by functions of $\LP\RR$ and $Z$ is $\tau$-continuous.

The paper is organized as follows. In Section~\ref{preliminaries}, we recall basic results and properties from convex geometry, the Monge--Ampère measure, the Legendre transform, and infimal convolution. Section~\ref{dual results} is devoted to the dual space of $\LP(\mathbb{R}^n)$ via the Legendre transform, where we derive the dual theorem corresponding to Theorem~\ref{equiv} and the dual corollary corresponding to Corollary~\ref{cor1}. Finally, in Section~\ref{Proof of the Theorem}, we present the proof of Theorem~\ref{equiv}, which relies on several intermediate results established there. The proofs of some technical auxiliary results are given in Appendix~\ref{apendice}.

\section{Background}\label{preliminaries}
In this section, we collect several tools and basic results that will be used throughout the paper, including definitions, properties of convex functions, and classical results.

\subsection{Preliminaries}\label{preli}
Following \cite{rockafellar2009variational},  a function $u\in\Conv\RR$ is called \textit{piecewise affine} if $\dom u$ can be written as the union of finitely many polyhedral sets, relative to each of which $u(x)$ is  affine. A set $C\subset\mathbb{R}^n$  is said to be  \textit{polyhedral}  if it can be expressed as the intersection of finitely many closed halfspaces. We denote by $\PA\RR$ the set of all piecewise affine functions.   A function $w\in \Conv\RR$ is called \textit{piecewise linear-quadratic} if $\dom w$ can be expressed  as the union of finitely many polyhedral sets, on each of which $w$ has the form
$$
w(x) = \frac{1}{2}\langle x,Ax\rangle + \langle y, x \rangle + c, \quad x\in \dom w,
$$
for some $c\in \mathbb{R}$, $y\in \mathbb{R}^n$, and a symmetric positive semidefinite matrix $A \in \mathbb{R}^{n \times n}$.

By \cite[Lemma 14]{baeta1d}, we obtain the following result.

\begin{lema}\label{piecewise}
For every $u\in \LC\RR$, there exists a sequence $u_k\in \PA(\mathbb{R}^n)\cap \LP\RR$, $k\in\mathbb{N}$,  that is $\tau$-convergent to $u$.
\end{lema}

One natural class of functions in $\LC\RR$ is formed by the \textit{indicator functions} of convex bodies $K\in\mathcal{K}^n$
\begin{align*}
\I_K(x) = 
\begin{cases}
0,& x \in K,\\
\infty,& x \notin K.
\end{cases}
\end{align*}

We say that $Z: \LP\RR\rightarrow \mathbb R$ is  \textit{simple} if $Z(u)=0$ for every  $u\in \LP\RR$ such that $\dom u$ is not full-dimensional. Let $\mathcal{P}^n$ denote the set of polytopes in $\mathbb{R}^n$. The following result from \cite[Lemma 13]{baeta1d} will be used in the proof of Theorem~\ref{equiv}.

\begin{lema}\label{non-positive}
Let $Z: \LP\RR\rightarrow \mathbb R$ be  a simple, $\tau$-upper semicontinuous, equi-affine  and dually epi-translation invariant valuation. Then there exists a constant $d\in\mathbb{R}$ such that 
\begin{align*}
   Z(\I_P) = d\, V_n(P)
\end{align*}
for every polytope $P\in \mathcal{P}^n$.
\end{lema}

Strictly speaking, Lemma~13 in~\cite{baeta1d} is stated for $\LC(\mathbb{R}^n)$. However, its proof only involves indicator functions of polytopes and does not use any additional properties of $\LC(\mathbb{R}^n)$. Therefore, it remains valid when restricted to $\LP(\mathbb{R}^n)$.

 We say that a valuation $\Phi: \mathcal{P}^n\to \mathbb{R}$ is \textit{simple} if $\Phi(P)=0$ for every polytope $P\subset \mathbb{R}^n$ of dimension at most $n-1$.  See \cite[Theorem 14.1]{mullen2} and the references therein for the following well-known result.

\begin{teo}\label{mullen}
Let $\Phi$ be a translation invariant, non-negative, simple valuation on $\mathcal{P}^n$. 
Then there exists a non-negative constant $d \in \mathbb{R}$ such that
$$
    \Phi(P) = d\, V_n(P)
$$
for every polytope $P \in \mathcal{P}^n$.
\end{teo}

To prove Proposition \ref{prop10.3} in Section \ref{Proof of the Theorem}, we will use the following version of the Vitali covering theorem.  A collection of sets $\mathcal{C}$ is called a \textit{Vitali covering} for $B \subset \mathbb{R}^n$ with respect to $V_n$ if, for each $x\in B$ and $\delta> 0$, there exists a $U\in \mathcal{C}$ such that $x\in U$ and $0<V_n(U)\leq \delta$. We say that  $\mathcal{C}$ is  a  \emph{regular family} for $B \subset \mathbb{R}^n$  with respect to $V_n$ if there exists a constant $d>0$ such that $\diam (U)^n\leq d \, V_n(U)$ for every set $U\in \mathcal{C}$,  where $\operatorname{diam}$ denotes the diameter.

\begin{teo}[\cite{falconer1985geometry}, Theorem 1.10]\label{VCT}
Let $ B\subset \mathbb{R}^n$ be a bounded set, and let $\mathcal{C}$ be a regular family of closed sets for $B$  with respect to $V_n$ that is a Vitali covering for $B$. Then, for any $\varepsilon > 0$, there exists a finite collection of pairwise disjoint sets $U_1, \dots, U_m \in \mathcal{C}$ such that
$$
V_n\Big(B \setminus \bigcup_{i=1}^m U_i \Big) < \varepsilon.
$$
\end{teo}

We  denote by  $\conv(A)$ the convex hull of a set $A\subset\mathbb{R}^n$, that is, the smallest convex set containing $A$.   For functions $f,g:\mathbb{R}^n \to (-\infty,\infty]$, we set
$$\{f\neq g\}=\{x\in\mathbb{R}^n\mid f(x)\neq g(x)\}\qquad \mbox{and} \qquad \{f=g\}=\{x\in\mathbb{R}^n\mid f(x)= g(x)\}.$$
A family $\mathcal{P}$ of pairwise disjoint non-empty subsets of $A \subset \mathbb{R}^n$ is called a \textit{partition} of $A$ if the union of all sets in $\mathcal{P}$ equals $A$. 

\medskip

The following theorem provides a concrete characterization of epi-convergence, showing its equivalence with a more familiar notion of convergence. 

\begin{teo}[\cite{rockafellar2009variational}, Theorem 7.17]\label{equivalence}
Let $u_k$ be a sequence in $\Conv\RR$, $k\in\mathbb{N}$,  that epi-converges to $u\in\Conv\RR$. If $\dom u$ has non-empty interior, the following are equivalent
\begin{itemize}
    \item[(1)] $u_k$ epi-converges to $u$;
    \item[(2)] $u_k$ converges uniformly to $u$ on every compact set $C$ that does not contain a boundary point of $\dom u$.
\end{itemize}
\end{teo}

\subsection{Monge--Ampère Measure}
For every $u\in\Conv\RR$, the \textit{subdifferential} of $u$ at $x\in \dom u$ is defined by
\begin{align*}
    \partial u(x)=\{y\in\mathbb{R}^n\mid u(z)\geq u(x)+\langle y, z-x\rangle \ \mbox{for all} \ z\in\mathbb{R}^n\},
\end{align*}
and we set $\partial u(x)=\emptyset$ for $x\not\in\dom u$.
For a convex function $u:\mathbb{R}^n\rightarrow (-\infty,\infty]$, differentiability at a point $x \in  \mathbb{R}^n$ is equivalent to  the subdifferential of $u$ at $x$ being a singleton.  Moreover,  $L$  is the Lipschitz constant of  $u\in \Conv\RR$ in the interior of its domain  if and only if $\|y\|\leq L$ for every $y\in\partial u(x)$, where $\|y\|=\sqrt{\langle y,y \rangle}$ is the usual Euclidean norm of $y\in\mathbb{R}^n$ and $x\in \interior (\dom u)$ (see, e.g., \cite[Theorem 3.61]{beck2017first}). 

By formula (2.4) in \cite{schneider2014convex}, the subdifferential of $\I_K$ at each boundary point $x$ of $K$ is the normal cone of $K$ at $x$, which is given by
$$
N_K(x) = \{ y \in \mathbb{R}^n \mid \langle y, z - x \rangle \le 0 \ \text{for all } z \in K \}.
$$
In particular, there are functions that are Lipschitz continuous in the interior of their domains but whose subdifferentials are unbounded at boundary points. For example, $\partial \I_{[0,1]}(1)=[0,\infty)$.

An important Radon measure associated with a  convex function $v:\mathbb{R}^n\rightarrow\mathbb{R}$ is the \textit{Monge--Ampère measure},  defined by
\begin{align*}
    \MA(v; B)= V_n(\partial v(B)), 
\end{align*}
for any Borel set $B \subseteq \mathbb{R}^n$, where
$$\partial v(B)= \bigcup_{x\in B}\partial v(x).$$ 
The notion of Monge--Ampère measure is fundamental in the definition of weak or generalized solutions of the Monge--Ampère  equation (see,  e.g., \cite{figalli2017monge, Trudinger-Wang-MA}).

\subsection{Legendre Transform}
For a convex function $u\in \Conv(\mathbb{R}^n)$, the function 
\begin{align*}
    u^*(x)=\sup_{y\in\mathbb{R}^n}  \big(\langle x,y \rangle - u(y)\big), \quad x \in\mathbb{R}^n,
\end{align*}
is called the \textit{Legendre transform} (or \textit{convex conjugate}) of $u$.  By \cite[Theorem 1.4.1]{balestro2024convexity}, we have $u^{**}=u$. Note that convex conjugation is continuous with respect to\ epi-convergence,
\begin{equation}\label{epi}
   u_k \text{ is epi-convergent to } u\,\,\Leftrightarrow\,\, u^*_k \text{ is epi-convergent to } u^*
\end{equation}
for $u, u_k\in \Conv(\mathbb{R}^n)$, $k\in\mathbb{N}$,  (see \cite[Theorem 11.43]{rockafellar2009variational}). 
\medskip

Let $u \in \Conv(\mathbb{R}^n)$, $c \in \mathbb{R}$, $\varphi \in \V\N$, and $y \in \mathbb{R}^n$. Denote by $\tau_{y}(x)=x+y$ the translation by $y$ and by $l_y(x)=\langle x,y\rangle$ the linear function associated with $y$. A straightforward computation yields
\begin{itemize}
    \item[(1)] $(u+c)^*(x)=u^*(x)-c$;
    \item[(2)] $(u\circ\varphi)^*(x)=(u^*\circ\varphi^{-t})(x)$;
    \item[(3)] $(u\circ\tau_y^{-1})^*(x) = u^*(x)+l_y(x)$;
    \item[(4)] $(u+l_y)^*(x)= u^*\circ\tau_y^{-1}(x)$.
\end{itemize}
Furthermore, if $u,v\in \Conv(\mathbb{R}^n)$ are such that $u\wedge v\in \Conv(\mathbb{R}^n)$, then  $u^*\wedge v^*\in \Conv(\mathbb{R}^n)$,  and
\begin{align}\label{duall}
(u \wedge v)^* = u^* \vee v^*, 
\qquad 
(u \vee v)^* = u^* \wedge v^*
\end{align}
(see, e.g., \cite[Lemma 3.4 and Proposition 3.5]{colesanti2020hessian}).

By \cite[Theorem 11.14]{rockafellar2009variational}, if $u\in \Conv(\mathbb{R}^n)$ has compact domain, then its Legendre transform $u^*$ is finite on $\mathbb{R}^n$, and $u$ is piecewise affine if and only if $u^*$ is piecewise affine.

\subsection{Infimal Convolution}\label{bauschke2011convex}
For any two functions $u,v\in \Conv(\mathbb{R}^n)$, the \textit{infimal convolution} $u\,\square\, v$ is defined by
\begin{align}\label{inf-}
   \left(u\,\square\, v\right)(x)=\inf_{y\in\mathbb{R}^n}\left(u(y)+v(x-y)\right), \quad x\in \mathbb{R}^n.
\end{align}
The epigraph of $u\in \Conv(\mathbb{R}^n)$ is the non-empty, closed, convex set
\begin{equation*}
    \epi(u)=\{(x, t)\in\mathbb{R}^{n+1}\mid  u(x)\leq t \}.
\end{equation*} 
Geometrically, the infimal convolution corresponds to the Minkowski sum of epigraphs,
$$
\epi(u\,\square\, v)= \epi(u)+\epi(v)
$$
where $A+B=\{a+b \mid a\in A,\ b\in B\}$.
In particular, if $v(x)=\frac{\lambda}{2}\|x\|^2$ with $\lambda>0$, then $u\,\square\, v$ is a smooth convex function (indeed, $C^1$ with Lipschitz gradient), even when $u$ is non-smooth and contains corners; see, e.g., \cite[Theorem~2.26]{rockafellar2009variational}.

By \cite[Proposition 2.1 and Proposition~2.6]{colesanti2013first}, for every $u,v\in\LP\RR$, we have 
$$
\dom(u\,\square\, v)= \dom u + \dom v, 
\quad u\,\square\, v\in\LP\RR,
$$
and the following relation holds
\begin{align}\label{inf-sum}
    (u\,\square\, v)^* = u^*+v^*.
\end{align}

Following \cite{ludwig1999characterization}, recall that a convex body in $\mathbb{R}^n$ is called a \textit{cylinder set} if it can be written as the Minkowski sum of an at most $(n-1)$-dimensional convex body $K$ and a closed line segment $J$. We define an analogous notion for convex functions.  
We say that a function $u\in\LC\RR$ is a \textit{cylinder function} if there exist $w \in \LC(\mathbb{R}^n)$, a closed interval $J\subset \mathbb{R}^n$, and an affine function $v$ such that
$$
u=w\,\square\, (v+\I_{J}),
$$
where $\dom w$ is a convex body of dimension at most $(n-1)$.

A particularly important example of infimal convolution is the \textit{Moreau envelope}, defined for $\lambda>0$ and $u\in \Conv(\mathbb{R}^n)$ by
\begin{align*}
    u_{\lambda}= u\,\square\, \frac{\lambda}{2}\|\cdot\|^2
\end{align*}
(see, e.g., \cite{bauschke2011convex,rockafellar2009variational}). By \cite[Proposition 12.15]{bauschke2011convex}, the infimum in \eqref{inf-} that defines $u_{\lambda}$ is uniquely attained. Moreover, by \cite[Theorem 1.25]{rockafellar2009variational},  $u_\lambda$ epi-converges to $u$ as $\lambda \nearrow \infty$, and by \cite[Theorem 2.26]{rockafellar2009variational}, the Moreau envelope $u_\lambda$  is convex and continuously differentiable.

We will work with a variant of the Moreau envelope that preserves both the compactness and the polytopal structure of the domain. Let $\lambda,\mu>0$, and denote  $C=[-1,1]^n$. 
For $u\in \LP\RR$, define
\begin{align}\label{ulm}
    u_{\lambda,\mu}= u\,\square\, \left(\frac{\lambda}{2}\|\cdot\|^2+\I_{\mu C}\right).
\end{align}
Note that $\dom u_{\lambda,\mu}=\dom u+\mu C$, which is a polytope. Since $u$ is Lipschitz on $\dom u$, a straightforward argument based on the definition of the infimal convolution shows that $u_{\lambda,\mu}$ is Lipschitz on $\dom u+\mu C$. Hence, $u_{\lambda,\mu}\in\LP(\mathbb{R}^n)$.
\medskip

For $x=(x_1,\dots,x_n)\in \mathbb{R}^n$ and each $i=1,\dots,n$, define
\begin{align}\label{q_i}
 q_i(x)=\frac{x_i^2}{2},   
\end{align}
and let $e_1,\dots,e_n$ be the canonical basis of $\mathbb{R}^n$. 

Let $q(x)= \frac{\|x\|^2}{2}$ for $x\in\mathbb{R}^n$. This quadratic function is the unique function satisfying $q^* = q$, see, for example, \cite[Example 11.11]{rockafellar2009variational}. Moreover,
\begin{align*}
 q(x) = \sum_{i=1}^n \frac{x_i^2}{2} =   \sum_{i=1}^n q_i(x), 
\end{align*}
where  $x=(x_1,\dots,x_n)$. Thus, by \eqref{inf-sum},
\begin{align*}
    q= q^*= \left(\sum_{i=1}^n q_i\right)^*= q_1^*\,\square\, \cdots \,\square\, q_n^*,
\end{align*}
and consequently,
\begin{align}\label{ccc}
u_{\lambda,\mu} = u \,\square\, \lambda\left(q+\I_{\mu C}\right)= u\,\square\, \lambda\left(q_1^*\,\square\, \cdots \,\square\, q_n^*+\I_{\mu C}\right).
\end{align}

We write $\mathrm{span}\{v\}=\{tv:t\in\mathbb{R}\}$.
We first establish the following structural property.

\begin{lema}\label{qi-star}
For each $i=1,\dots,n$, one has
\[
q_i^*(x)=\frac{x_i^2}{2}+\I_{\mathrm{span}\{e_i\}}(x).
\]
\end{lema}

\begin{proof}
Fix $i\in\{1,\dots,n\}$. By definition,
\[
q_i^*(x)=\sup_{y\in\mathbb{R}^n}\bigl(\langle x,y\rangle-\tfrac{1}{2}y_i^2\bigr).
\]
If $x\notin \mathrm{span}\{e_i\}$, then the supremum is $\infty$. Hence $q_i^*(x)<\infty$ only if $x\in \mathrm{span}\{e_i\}$. In this case, writing $x=x_i e_i$, we obtain
\[
q_i^*(x)=\sup_{y_i\in\mathbb{R}}\left(x_i y_i-\tfrac{1}{2}y_i^2\right)=\frac{x_i^2}{2}.
\]
This yields the claimed formula.
\end{proof}

\begin{lema}\label{decomposition}
Let $u\in\LP(\mathbb{R}^n)$ and let $\lambda,\mu>0$. Then
\begin{align*}
u_{\lambda,\mu}
= u\,\square\, \lambda \left[
\left(q_1+\I_{[-\mu,\mu]\times \{0\}^{n-1}} \right)
\,\square\,
\left(q_2+\I_{\{0\}\times [-\mu,\mu]\times \{0\}^{n-2}} \right)
\,\square\, \cdots \,\square\,
\left(q_n+\I_{\{0\}^{n-1}\times [-\mu,\mu]} \right)
\right].
\end{align*}
\end{lema}

\begin{proof}
For $i=1,\dots,n$, set
\(
S_i=\{te_i:t\in[-\mu,\mu]\}.
\)
By Lemma~\ref{qi-star}, the  domain of $q_i^*$ is the $i$-th coordinate axis. Thus, in order to obtain the restriction to the cube $\mu C$, it is enough to restrict each coordinate component to the segment $S_i$.

We first show that
\[
q+\I_{\mu C}
=
(q_1+\I_{S_1})\,\square\,\cdots\,\square\,(q_n+\I_{S_n}).
\]
Taking Legendre transforms and using \eqref{inf-sum}, the right-hand side becomes
\[
\sum_{i=1}^n (q_i+\I_{S_i})^*(x)
=
\sum_{i=1}^n \sup_{t_i\in[-\mu,\mu]}
\left(t_i x_i-\frac{t_i^2}{2}\right).
\]
On the other hand,
\[
(q+\I_{\mu C})^*(x)
=
\sup_{y\in\mu C}
\left(\langle x,y\rangle-\frac{\|y\|^2}{2}\right)
=
\sum_{i=1}^n \sup_{y_i\in[-\mu,\mu]}
\left(x_i y_i-\frac{y_i^2}{2}\right).
\]
Thus the Legendre transforms coincide. Since all functions involved are proper lower semicontinuous convex functions, taking Legendre transforms again yields
\[
q+\I_{\mu C}
=
(q_1+\I_{S_1})\,\square\,\cdots\,\square\,(q_n+\I_{S_n}).
\]
Using the definition
\[
u_{\lambda,\mu}
=
u\,\square\,\lambda(q+\I_{\mu C}),
\]
the result follows.
\end{proof}

A classical construction associates to each convex body $K \in \mathcal{K}^{n+1}$ a convex function with compact domain in $\mathbb{R}^n$ via
\begin{align*}
 \lfloor K \rfloor(x)=\inf_{(x,t)\in K}t   
\end{align*}
(see, for example,  \cite{knoerr2024valuations}). Conversely, for a convex function $u \in \LP(\mathbb{R}^n)$, one can associate the convex body
\begin{align}\label{set_func}
    K^u=\epi(u-M_u)\cap R_{\mathbb{R}^n}(\epi(u-M_u))+M_ue_{n+1}
\end{align}
where $M_u=\max_{x\in \dom u}u(x)$ and $R_{\mathbb{R}^n}$ denotes the reflection through the hyperplane $\mathbb{R}^n = e_{n+1}^\perp$. Moreover,
$$\lfloor K^u \rfloor = u.$$

\begin{lema}\label{inf_sum}
Let $u_1,u_2 \in \LC(\mathbb{R}^n)$. Then
\[
K^{u_1 \square u_2}= K^{u_1}+K^{u_2}.
\]
\end{lema}

\begin{proof}
Let $c_i=\max_{x\in \dom u_i} u_i(x)$ for $i=1,2$. By \cite[Lemma 5.1]{mussnig2026inequalities},
\[
K^{(u_1-c_1)\square (u_2-c_2)} = K^{u_1-c_1}+K^{u_2-c_2}.
\]
By definition of infimal convolution,
\[
(u_1-c_1)\square(u_2-c_2) = (u_1\square u_2) - (c_1+c_2),
\]
and by \eqref{set_func},
\[
K^{u_i-c_i} = K^{u_i} - c_i e_{n+1}, \quad i=1,2.
\]
Therefore,
\begin{align*}
K^{u_1\square u_2}
&= K^{(u_1-c_1)\square (u_2-c_2)} + (c_1+c_2)e_{n+1} \\
&= K^{u_1-c_1} + K^{u_2-c_2} + (c_1+c_2)e_{n+1} \\
&= K^{u_1}+K^{u_2}.
\end{align*}
This concludes the proof.
\end{proof}

A set $S\subset\mathbb{R}^{n+1}$ is called a \textit{zonotope} if $S= J_1+\cdots+J_m$ for some line segments $J_1,\dots,J_m$, with respect to Minkowski sum. If a convex body $T\in\mathcal{K}^{n+1}$ can be approximated in the Hausdorff metric by zonotopes, then $T$ is called a \textit{zonoid}. We say that a set $A$ is centrally symmetric if there exists a point $x_0 \in A$ such that
\(
A - x_0 = -(A - x_0).
\) By \cite[Corollary 3.5.7]{schneider2014convex}, every centrally symmetric convex body in $\mathbb{R}^2$ is a zonoid. Moreover, by \cite[Theorem 3.5.2]{schneider2014convex}, a convex polytope is a zonotope if and only if all its $2$-dimensional faces are centrally symmetric.

\begin{lema}\label{lemma:zonotope-approx}
Let $\mu>0$ and define $v:\mathbb{R}\to(-\infty, \infty]$ by
\[
v(x)=\frac{x^2}{2}+\I_{[-\mu,\mu]}(x).
\]
Let $K^v$ be the convex body associated with $v$ via \eqref{set_func}. Then there exists a sequence of zonotopes $P_k \subset \mathbb{R}^2$ such that $P_k \to K^v$ in the Hausdorff metric as $k\to\infty$, and such that, if
\(
w_k=\lfloor P_k\rfloor,
\)
then $w_k\in\LP(\mathbb{R})$ and $w_k\seq v$.
\end{lema}

\begin{proof}
By \eqref{set_func}, we have
\[
K^v=
\left\{
(x,t)\in\mathbb{R}^2:
x\in[-\mu,\mu],\;
\frac{x^2}{2}\le t\le \mu^2-\frac{x^2}{2}
\right\}.
\]
In particular, $K^v$ is a centrally symmetric convex body in $\mathbb{R}^2$.

For each $k\in\mathbb{N}$, let
\[
x_i^k=-\mu+\frac{2\mu i}{k},
\qquad i=0,\dots,k,
\]
and let $w_k:\mathbb{R}\to(-\infty,\infty]$ be the function that agrees on
$[-\mu,\mu]$ with the piecewise affine interpolation of $x^2/2$ at the points
$x_i^k$, $i=0,\dots,k$, and is equal to $\infty$ outside
$[-\mu,\mu]$. Then $w_k\in\LP(\mathbb{R})$ for every $k\in\mathbb{N}$. 

Define
\[
P_k=K^{w_k}.
\]
Equivalently,
\[
P_k=
\left\{
(x,t)\in\mathbb{R}^2:
x\in[-\mu,\mu],\;
w_k(x)\le t\le \mu^2-w_k(x)
\right\}.
\]
Then $P_k$ is a convex polygon which is centrally symmetric with respect to
$(0,\mu^2/2)$ for every $k\in\mathbb{N}$. Hence $P_k$ is a zonotope in
$\mathbb{R}^2$ for every $k\in\mathbb{N}$, since every centrally symmetric convex polygon in $\mathbb{R}^2$ is a zonotope. Moreover, by construction,
\[
\lfloor P_k\rfloor=w_k
\]
for every $k\in\mathbb{N}$.

Since $w_k$ converges uniformly to $v$ on $[-\mu,\mu]$ as $k\to\infty$, and $w_k$ and $v$ have the same domain for every $k\in\mathbb{N}$, it follows that
\[
P_k\to K^v
\]
in the Hausdorff metric as $k\to\infty$.

It remains to verify that the Lipschitz constants of $w_k$, $k\in\mathbb{N}$, are uniformly bounded. For every $k\in\mathbb{N}$ and every $i=0,\dots,k-1$, the slope of $w_k$ on the interval $[x_i^k,x_{i+1}^k]$ is
\[
\frac{w_k(x_{i+1}^k)-w_k(x_i^k)}
{x_{i+1}^k-x_i^k}
=
\frac{(x_{i+1}^k)^2-(x_i^k)^2}
{2(x_{i+1}^k-x_i^k)}
=
\frac{x_{i+1}^k+x_i^k}{2}.
\]
Since $x_i^k\in[-\mu,\mu]$ for every $i=0,\dots,k$, all these slopes belong to
$[-\mu,\mu]$. Therefore, the Lipschitz constants of $w_k$, $k\in\mathbb{N}$, are bounded by $\mu$. Hence $w_k\seq v$, completing the proof.
\end{proof}

By Lemma~\ref{lemma:zonotope-approx}, we obtain the following.

\begin{lema}\label{1-dim}
For each $i=1,\dots,n$ and $x=(x_1,\dots,x_n)\in \mathbb{R}^n$, let
\[
v_i(x)=\frac{x_i^2}{2}+\I_{\{0\}^{i-1}\times[-\mu,\mu]\times\{0\}^{n-i}}(x).
\]
Then there exist line segments
$J_1^{k,i},\dots,J_{m_k^{(i)}}^{k,i} \subset \mathbb{R}^2$
such that
\[
w_k^{(i)} = \lfloor J_1^{k,i} \rfloor \,\square \cdots \square\, \lfloor J_{m_k^{(i)}}^{k,i} \rfloor,
\quad \text{and} \quad
w_k^{(i)} \seq v_i \ \text{as } k\to\infty.
\]
Moreover, the Lipschitz constants of $w_k^{(i)}$ are uniformly bounded in $k$ and $i$.  Finally, for $u\in \LP\RR$, defining
\(
\tilde{w}_k = w_k^{(1)} \square \cdots \square w_k^{(n)},
\)
we have
\[
u \square \lambda \tilde{w}_k \seq u_{\lambda,\mu},
\]
where $u_{\lambda,\mu}$ is given by \eqref{ulm}.
\end{lema}

\begin{proof}
Observe that each $v_i$ is obtained from $v_1$ by a permutation of coordinates, hence all cases are equivalent up to an orthogonal change of variables. By Lemma~\ref{lemma:zonotope-approx}, for each $k\in\mathbb{N}$ there exists a zonotope $P_k^{(i)} \subset \mathbb{R}^{n+1}$ such that the associated functions $w_k^{(i)}=\lfloor P_k^{(i)} \rfloor$ satisfy
\[
w_k^{(i)} \seq v_i \quad \text{as } k\to\infty,
\]
and their Lipschitz constants are uniformly bounded in $k$ and $i$. Since each $P_k^{(i)}$ is a zonotope, there exist line segments $J_1^{k,i},\dots,J_{m_k^{(i)}}^{k,i}\subset\mathbb{R}^2$ such that
\[
P_k^{(i)} = J_1^{k,i}+\cdots+J_{m_k^{(i)}}^{k,i}.
\]
Denoting by $\lfloor J_j^{k,i}\rfloor$ the associated convex functions, Lemma~\ref{inf_sum} yields
\[
w_k^{(i)} = \lfloor J_1^{k,i} \rfloor \,\square \cdots \square\, \lfloor J_{m_k^{(i)}}^{k,i} \rfloor,
\]
and $w_k^{(i)} \seq v_i$ as $k\to\infty$.

Define
\[
\tilde{w}_k = w_k^{(1)}\,  \square \cdots \square \, w_k^{(n)}\in \LP(\mathbb{R}^n).
\]
Using this approximation together with Lemma~\ref{decomposition}, we obtain
\begin{align*}
u\,\square\, \lambda \tilde{w}_k \seq u_{\lambda,\mu} \quad \text{as } k\to \infty.
\end{align*}
This completes the proof of the lemma.
\end{proof}

Since $\frac{\lambda}{2}\|\cdot\|^2+\I_{\mu C}$ epi-converges to $\I_{\{0\}}$ as $\lambda\rightarrow \infty$, it follows from \eqref{epi} that  $\left(\frac{\lambda}{2}\|\cdot\|^2+\I_{\mu C}\right)^*$ epi-converges to the constant function $l\equiv 0$. By \eqref{inf-sum},  we have
\begin{align*}
    (u_{\lambda,\mu})^*= \left(u\,\square\, \left(\frac{\lambda}{2}\|\cdot\|^2+\I_{\mu C}\right)\right)^*= \left(u^*+\left(\frac{\lambda}{2}\|\cdot\|^2+\I_{\mu C}\right)^*\right).
\end{align*}
Consequently, $(u_{\lambda,\mu})^*$ epi-converges to $u^*$, which is equivalent to $u_{\lambda,\mu}$ epi-converging to $u$ as $\lambda\rightarrow \infty$ with $\mu>0$. However, $u_{\lambda,\mu}$ is  not $\tau$-convergent to $u$ as $\lambda\rightarrow \infty$, since the Lipschitz constants are not uniformly bounded. To address this, instead of letting $\lambda \to \infty$, we fix $\lambda>0$ and let $\mu\to 0$, as stated in the following result.

\begin{lema}
Let $u\in \LP(\mathbb{R}^n)$ and fix $\lambda>0$. Then $u_{\lambda,\mu}$ is $\tau$-convergent to $u$ as $\mu \to 0$.
\end{lema}

\begin{proof}
Since $\I_{\mu C}$ epi-converges to $\I_{\{0\}}$ as $\mu \to 0$, it follows that
\(
\frac{\lambda}{2}\|\cdot\|^2 + \I_{\mu C}
\)
epi-converges to $\I_{\{0\}}$.
By \eqref{inf-sum}, we obtain that $u_{\lambda,\mu}$ epi-converges to
\(u\). Finally, since the functions $u_{\lambda,\mu}$ have uniformly bounded Lipschitz constants for fixed $\lambda$, epi-convergence implies $\tau$-convergence, concluding the proof.
\end{proof}

\begin{lema}\label{gradient}
Let $u\in \LP(\mathbb{R}^n)$ and let $\lambda,\mu>0$. 
Then, for every $x_0\in\dom u_{\lambda,\mu}$, where $u_{\lambda,\mu}$ is defined by \eqref{ulm}, there exists $y_0 \in \dom u$ at which the infimum in the definition of $u_{\lambda,\mu}(x_0)$ is attained, and there exists a quadratic function 
$q_{x_0}:\mathbb{R}^n\to\mathbb{R}$ defined by

\[
q_{x_0}(x) = u(y_0) + \frac{\lambda}{2}\|x - y_0\|^2
\]
such that $q_{x_0}(x_0)=u_{\lambda,\mu}(x_0)$ and $q_{x_0}(x)\ge u_{\lambda,\mu}(x)$ for all $x\in y_0+\mu C\subset \dom u_{\lambda,\mu}$. In particular, $x_0\in y_0+\mu C$.
\end{lema}

\begin{proof}
By definition,
$$
u_{\lambda,\mu}(x)
=\inf_{\substack{y\in\dom u \\ x-y\in\mu C}}\Bigl(u(y)+\frac{\lambda}{2}\|x-y\|^2\Bigr), \qquad x\in\mathbb{R}^n.
$$
Since $x_0\in\dom u_{\lambda,\mu}$, the  set
$$
S_{x_0}=\{y\in\dom u\mid \;x_0-y\in\mu C\}
$$
is non-empty and compact, and the function 
$y\mapsto u(y)+\frac{\lambda}{2}\|x_0-y\|^2$ 
is lower semicontinuous. Hence the infimum is attained, i.e., there exists $y_0\in S_{x_0}$ such that
$$
u_{\lambda,\mu}(x_0)
= u(y_0)+\frac{\lambda}{2}\|x_0-y_0\|^2.
$$

Define the quadratic function
$$
q_{x_0}(x)=u(y_0)+\frac{\lambda}{2}\|x-y_0\|^2, \qquad x\in\mathbb{R}^n.
$$
Then $q_{x_0}(x_0)=u_{\lambda,\mu}(x_0)$. Moreover, for any $x$ such that $x-y_0\in\mu C$, that is, $x\in y_0+\mu C$,
we have
$$
u_{\lambda,\mu}(x)
=\inf_{\substack{y\in\dom u \\ x-y\in\mu C}}\Bigl(u(y)+\tfrac{\lambda}{2}\|x-y\|^2\Bigr)
\le u(y_0)+\frac{\lambda}{2}\|x-y_0\|^2
=q_{x_0}(x).
$$
Therefore,
$$
q_{x_0}(x)\ge u_{\lambda,\mu}(x)
$$
for all $x\in y_0+\mu C$. Moreover, since $x_0 - y_0 \in \mu C$, we have $x_0 \in y_0 + \mu C$.    This completes the proof. 
\end{proof}

\section{Dual Results}\label{dual results}
Let $K\subset\mathbb{R}^{n+1}$ be a non-empty closed convex set with support function
$h_K:\mathbb{R}^{n+1}\to (-\infty, \infty]$ defined by
$$h_{K}(x,t)=\sup_{(y,y_{n+1})\in K}\langle (x, t), (y,y_{n+1})\rangle$$
for $x\in \mathbb{R}^{n}$ and $t\in\mathbb{R}$. Note that  for any $u\in \Conv(\mathbb{R}^n)$, we have
\begin{equation}\label{minus}
    u^*(y)=\sup_{x\in\mathbb{R}^n}\langle (y, -1),(x,u(x))\rangle = h_{\epi(u)}(y,-1).
\end{equation}
If $u\in \Conv(\mathbb{R}^n)$ has compact domain and $K = K^u$ as in \eqref{set_func},  we obtain
\begin{equation}\label{minus*}
    h_{K}(\cdot,-1) = h_{\epi(u)}(\cdot,-1)\in \Conv\finite.
\end{equation}
In particular,  if $\dom u$ is a polytope, then the projection of $K$ onto $\mathbb{R}^n=e_{n+1}^{\perp}$ is also a polytope.

Define
\begin{align*}
    \Conv_{\MAh}\finite=
    \{h_{K}(\cdot, -1)\in\Conv\finite\mid \, K\in\mathcal{K}^{n+1},\, \supp(\MA(h_{K}(\cdot, -1);\cdot)) \text{ is compact}\},
\end{align*}
where supp denotes the support of a measure. We also consider the subset
\begin{align*}
    \Conv_{\MAp}\finite
    =
    \{h_{K}(\cdot, -1)\in \Conv_{\MAh}\finite\mid\,  \text{ the projection of $K$ onto $\mathbb{R}^n$ is a polytope}\}.
\end{align*}
We say that a sequence $v_k\in \Conv_{\MAh}\finite$, $k\in\mathbb{N}$,  is $\tau^*$-convergent to $v\in \Conv_{\MAh}\finite$ if the following conditions hold
\begin{enumerate}
\item[(i)] $v_k$ epi-converges to $v$;
\item[(ii)] there exists a compact set $C\subset\mathbb{R}^n$ such that $\supp(\MA(v_k;\cdot)),\supp(\MA(v;\cdot))\subseteq C$ for every $k$.
\end{enumerate}
(cf. \cite{baetaludwigsemicontinuity}). This definition also applies to $\Conv_{\MAp}\finite$.  Recall that by Theorem~\ref{equivalence}, for finite-valued convex functions, epi-convergence is equivalent to uniform convergence on compact sets. The terminology $\tau^*$-convergence is motivated by the duality between the sets $\Conv_{\MAh}\finite$ and $\LC\RR$, established via the Legendre transform, as follows.

\begin{lema}[\cite{baetaludwigsemicontinuity}, Lemma 3]\label{duality}
A sequence $u_k$ in $\LC\RR$, $k\in\mathbb{N}$, is $\tau$-convergent to $u \in \LC\RR$ if and only if $u_k^*$ and $u^*$ belong to 
$\Conv_{\MAh}\finite$ and $u_k^*$ is $\tau^*$-convergent to  $u^*$.
\end{lema}

By Lemma~\ref{duality}, together with \eqref{minus} and \eqref{minus*}, we  conclude that a sequence $u_k$ 
in $\LP\RR$, $k\in\mathbb{N}$,
is $\tau$-convergent to $u \in \LP\RR$ if and only if $u_k^*$ and $u^*$ belong to 
$\Conv_{\MAp}\finite$ and $u_k^*$ is $\tau^*$-convergent to  $u^*$.

Using \cite[Corollary 23.5.1]{rockafellar1997convex}, we have  
\begin{align}\label{sub}
    p\in \partial v(x) \quad  \Leftrightarrow \quad x \in \partial v^*(p),
\end{align}
for $v\in \Conv\finite$.

Let  $v\in \Conv\finite$ and $p\in\dom v^*$. By \cite[Lemma A.22]{figalli2017monge}, the set $\partial v(C)$ is compact for every  compact set $C\subset \mathbb{R}^n$. If $p\notin \partial v(\mathbb{R}^n)$, then $p\not\in\partial v(x)$ for every $x\in\mathbb{R}^n$. By \eqref{sub},  this implies 
\begin{align*}
    x\not\in \partial v^*(p) \quad \mbox{ for all }  x\in\mathbb{R}^n,
\end{align*}
that is, $\partial v^*(p)$ is empty. Hence $p$ lies on the boundary of $\dom v^*$. Therefore, 
\begin{align}\label{volum}
    V_n(\dom v^*)= V_n(\partial v(\mathbb{R}^n))+V_n(\dom v^*\setminus \partial v(\mathbb{R}^n)) = V_n(\partial v(\mathbb{R}^n))=\MA(v;\mathbb{R}^n)
\end{align}
for every $v\in \Conv_{\MAh}\finite$.

In \cite[Lemma 9]{baeta1d}, it is shown that the map $v^*\mapsto V_n(\dom v^*)$ is $\tau$-continuous, i.e., if $v_k^*\in \LC(\mathbb{R}^n)$, $k\in\mathbb{N}$, is $\tau$-convergent to $v^*\in \LC(\mathbb{R}^n)$, then $V_n(\dom v_k^*)\to V_n(\dom v^*)$ as $k\to\infty$. We note that this also follows from \eqref{volum}.

Finally, let $v\in\Conv_{\MAh}\finite$ and $\zeta \in\conc$. Then, as shown for example in \cite[equation~(11)]{baetaludwigsemicontinuity}, we have
\begin{align}\label{conjugate}
  \int_{\mathbb{R}^n} \zeta (\det \D v(x))\dif x  =\int_{\dom v^*} \tilde{\zeta} (\det \D v^*(x))\dif x,   
\end{align}
where $\tilde{\zeta}= \zeta(1/t)\,t$ for $t>0$, and $\tilde{\zeta}\in\conc$. 
\medskip

We now conclude, using \eqref{duall}, the relations derived above from it, and Lemma~\ref{duality},  that the functional  $Z^*: \Conv_{\MAp}\finite\rightarrow\mathbb{R}$ given by $Z^*(v)= Z(v^*)$ is a $\tau^*$-upper semicontinuous, equi-affine   and dually epi-translation invariant valuation if and only if $Z: \LP\RR\rightarrow\mathbb{R}$  is a $\tau$-upper semicontinuous, equi-affine    and dually epi-translation invariant valuation.  Combining this with \eqref{volum} and  \eqref{conjugate},  we deduce the following result from Theorem \ref{equiv}.

\begin{teo}\label{maintheorem}
A functional  $Z:\Conv_{\MAp}\finite\rightarrow\mathbb{R}$ is  a $\tau^*$-upper semicontinuous, equi-affine  and dually epi-translation invariant valuation  if and only if there are constants $c_0,c_1\in\mathbb{R}$ and a  function   $\zeta\in\conc$ such that
\begin{align*}
Z(v)= c_0+c_1\MA(v;\mathbb{R}^n)+\int_{\mathbb{R}^n}\zeta(\det \D v(x))\dif x
\end{align*}
for every $v\in \Conv_{\MAp}\finite$.
\end{teo}

We obtain the following result from Theorem~\ref{maintheorem}, or equivalently, by duality from Corollary~\ref{cor1}.

\begin{cor}
A functional  $Z:\Conv_{\MAh}\finite\rightarrow\mathbb{R}$ is a $\tau^*$-continuous, equi-affine and dually epi-translation  invariant valuation if and only if there are constants $c_0,c_1\in\mathbb{R}$ such that
$$Z(v)=c_0+c_1\MA(v;\mathbb{R}^n),$$ 
for every $v\in \Conv_{\MAh}\finite$.
\end{cor}

Since Theorem \ref{equiv} and Theorem \ref{maintheorem} are equivalent, it suffices to prove  Theorem \ref{equiv}.

\section{Proof of Theorem~\ref{equiv}}\label{Proof of the Theorem}

By the discussion following Theorem~\ref{equiv},  it is sufficient   to prove the necessity part of the theorem, namely, that if $Z: \LP\RR\rightarrow \mathbb{R}$ is a $\tau$-upper semicontinuous,  equi-affine  and dually epi-translation invariant valuation, then there exist constants $c_0,c_1$ and a function $\zeta\in\conc$ such that $Z$ is given by \eqref{maintheo}.  We will adapt elements from the proof of Theorem~\ref{cg} in \cite{ludwig1999characterization}. 

We begin by considering the case where
$Z:\LP\RR\rightarrow\mathbb{R}$ is a simple,
$\tau$-upper semicontinuous, equi-affine and dually epi-translation invariant valuation that vanishes on indicator functions of polytopes and on cylinder functions.

The following lemma establishes a basic non-negativity property.

\begin{lema}\label{non-negative}
Let $Z:\LP(\mathbb{R}^n)\to\mathbb{R}$ be a simple, $\tau$-upper semicontinuous, dually epi-translation invariant valuation that vanishes on indicator functions of polytopes. Then
\[
Z(u)\ge0
\]
for every $u\in\LP(\mathbb{R}^n)$.
\end{lema}

\begin{proof}
Using Lemma~\ref{piecewise}, there exists a sequence
$u_k\in\PA(\mathbb{R}^n)\cap\LP(\mathbb{R}^n)$, $k\in\mathbb{N}$, such that  $u_k$ is $\tau$-convergent to $u$. Since $Z$ vanishes on indicator functions of polytopes, is dually epi-translation invariant, simple, and $\tau$-upper semicontinuous, it follows that
\[
Z(u)\ge \limsup_{k\to \infty} Z(w_k)=0
\]
for every $u\in\LP(\mathbb{R}^n)$.
\end{proof}

Throughout the text, $q$ denotes the quadratic function defined by
\begin{align}\label{q}
    q(x)=\frac{1}{2}\sum_{i=1}^nx_i^2, \qquad x\in\mathbb{R}^n,
\end{align}  
and $C = [-1,1]^n$ denotes the unit cube centered at the origin with side length~$2$.

We now show that if $Z$ satisfies the assumptions above, then, for quadratic functions, the ratio $\frac{Z(aq+\I_P)}{V_n(P)}$ is independent of the choice of the polytope $P$, whenever $V_n(P)>0$.

\begin{lema}
Let $Z:\LP\RR\rightarrow\mathbb{R}$ be a simple $\tau$-upper semicontinuous, equi-affine and  dually epi-translation   invariant valuation that   vanishes on  indicator functions of polytopes. Then, for every $a>0$, 
\begin{align}\label{mudanca_dominio}
 \dfrac{Z(aq+\I_{P_1})}{V_n(P_1)}= \dfrac{Z(aq+\I_{P_2})}{V_n(P_2)}   
\end{align}
for all polytopes $P_1,P_2\subset\mathbb{R}^n$ with $V_n(P_1),V_n(P_2)>0$.
\end{lema}

\begin{proof}
Given $a>0$, define  the function $\Phi_a: \mathcal{P}^n\to \mathbb{R}$  by
$$\Phi_a(P)=Z(aq+\I_P),$$
where  $q(x)=\frac{1}{2}\sum_{i=1}^nx_i^2$, $x\in\mathbb{R}^n$.

By  Lemma \ref{non-negative}, the functional $\Phi_a$ is  non-negative. Moreover, since $Z$ is a simple valuation on $\LP\RR$, it follows that $\Phi_a$ is a simple valuation on $\mathcal{P}^n$. Let $y\in\mathbb{R}^n$ and consider the translation function $\tau_y(x)=x+y$. Note that
\begin{align*}
    q\circ\tau_y^{-1}(x)= q(x)+ \left(\tfrac{1}{2}\|y\|^2-\langle x,y\rangle\right),
\end{align*}
 for every $x\in\mathbb{R}^n$, and that $w_y(x)= \tfrac{1}{2}\|y\|^2-\langle x,y\rangle$ is an affine function. Since $Z$ is translation and  dually epi-translation invariant, we obtain
\begin{align*}
    Z(aq + \I_P)
      &= Z((aq + \I_P) \circ \tau_y^{-1}) \\
      &= Z(aq + aw_y + \I_{P+y}) \\
      &= Z(aq + \I_{P+y}).
\end{align*}
Hence, for every  $y\in\mathbb{R}^n$,
\begin{align*}
\Phi_a(P+y)  = Z(aq+\I_{P+y})=   Z(aq+\I_P)= \Phi_a(P). 
\end{align*}
This shows that $\Phi_a$ is translation invariant.

By Theorem \ref{mullen},  there exists a non-negative constant $d=d(a)\in\mathbb{R}$ such that 
$$\Phi_a(P)= d\,V_n(P),$$
and consequently,
\begin{align}\label{789}
 Z(aq+ \I_P)= d\,V_n(P)   
\end{align}
for every polytope $P\in\mathcal{P}^n$. In particular, \eqref{789} implies
\begin{align*}
 \dfrac{Z(aq+\I_{P_1})}{V_n(P_1)}= \dfrac{Z(aq+\I_{P_2})}{V_n(P_2)}   
\end{align*}
for all polytopes $P_1,P_2\subset\mathbb{R}^n$ with positive volume.
\end{proof}

Let $Z:\LP\RR\to\mathbb{R}$ be a simple, $\tau$-upper semicontinuous, equi-affine and dually epi-translation invariant valuation that vanishes on indicator functions of polytopes. Define the function $\zeta:[0,\infty)\rightarrow \mathbb{R}$ by
\begin{align}\label{zeta}
    \zeta (a)= \dfrac{Z(aq+\I_{C})}{V_n(C)}, 
\end{align}
where $a\geq 0$ and  $C=[-1,1]^n$. By \eqref{mudanca_dominio}, we then have
\begin{align*}
   Z(aq+\I_P) &=  \zeta(a)V_n(P),
\end{align*}
for every polytope $P\subset\mathbb{R}^n$ and every $a\geq 0$.  The following lemma generalizes this identity to arbitrary quadratic functions.

\begin{lema}
Let $Z:\LP(\mathbb{R}^n)\to\mathbb{R}$ be a simple $\tau$-upper semicontinuous, equi-affine and dually epi-translation invariant valuation that vanishes on indicator functions of polytopes. Let
$$
\tilde{q}(x)= \frac{1}{2}\langle x, Ax \rangle + l(x) + c,
$$
where $A$ is symmetric positive semidefinite with $\det A=a\geq 0$, $l:\mathbb{R}^n\to\mathbb{R}$ is linear, and $c\in\mathbb{R}$. Then
\begin{align}\label{mudanca_zeta}
   Z(\tilde{q}+\I_P)=\zeta(a)\,V_n(P)
\end{align}
for every polytope $P\subset\mathbb{R}^n$.
\end{lema}

\begin{proof}
The case $\tilde{q}=aq$ follows directly from the definition of $\zeta$.  For a general quadratic function $\tilde{q}$, there exists a linear transformation $T$ with $\det T=1$ such that
$$
\tilde{q}(x)= aq(Tx) + l(x) + c,
$$
up to an affine function.  Since $Z$ is equi-affine and dually epi-translation invariant, it is invariant under volume-preserving linear transformations and the addition of affine functions. Therefore,
\begin{align*}
Z(\tilde{q}+\I_P)
= Z(aq + \I_{T(P)}).
\end{align*}
By the definition of $\zeta$, this yields
$$
Z(\tilde{q}+\I_P)=\zeta(a)\,V_n(P),
$$
which completes the proof.
\end{proof}

By Lemma \ref{non-negative}, we then obtain the following result.

\begin{lema}\label{positive}
The function $\zeta$ given by \eqref{zeta} is non-negative.
\end{lema}

We will use  $\diag(a_1, \dots, a_n)$  to denote the $n \times n$ diagonal matrix with entries 
$$
x_{ij} = 
\begin{cases} 
a_i, & \text{if } i=j,\\
0, & \text{if } i\neq j,
\end{cases} 
\quad i,j=1,\dots,n.
$$

The properties of $Z$ as in Proposition~\ref{prop} allow us to conclude that $\zeta$ is concave. In dimension $n=1$, taking $q(x)=x^2/2$ and $P$ a closed interval, the function $\zeta$ is characterized by
\begin{align*}
    Z(aq+\I_P)=\zeta(a)\,V_1(P).
\end{align*}
The concavity of $\zeta$ in the one-dimensional case follows from \cite[Lemma 17]{baeta1d}.

\begin{lema}\label{conccave}
The function $\zeta$ given by \eqref{zeta} is concave.
\end{lema}

\begin{proof}
We begin with the case $n=2$. Let $0\leq s<a<r$ and $m\in \mathbb{N}$. Consider the points
$$p_i=\left(t_1, -t_2+\left(\frac{2t_2}{m}\right)i\right), \qquad i=0,\dots, m,$$
for $t_1,t_2>0$. 
Note that
$$p_0=(t_1,-t_2), \quad  p_m=(t_1,t_2),$$
and $p_i\in [(t_1,-t_2), (t_1,t_2)]$ for every $i=1, \dots, m-1$. In the two-dimensional case described above, let  $R=[-t_1,t_1]\times[-t_2,t_2]$ and  
$$q^a(x_1,x_2)= ax_1^2+x_2^2.$$
We  approximate the piecewise linear-quadratic function $q^a+\I_R$ by  suitable piecewise linear-quadratic functions as follows. Consider the functions
\begin{align*}
q_i^s(x_1,x_2)= ax_1^2 +\frac{s}{a}x_2^2 + 2\left(1-\frac{s}{a}\right)\left(-t_2+\left(\frac{2t_2}{m}\right)i\right)x_2-\left(1-\frac{s}{a}\right)\left(-t_2+\left(\frac{2t_2}{m}\right)i\right)^2,
\end{align*}
for every $i=0,\dots, m-1$. Note that
$$q_i^s(p_i)=q^a(p_i), \, \nabla q_i^s(p_i)= \nabla q^a(p_i) \ \mbox{ and } (q_i^s+\I_R)(x_1,x_2)  \leq (q^a+\I_R)(x_1,x_2) \mbox{ for all } (x_1,x_2)\in \mathbb{R}^2.$$

Next, we look for functions of the form
$$q_i^r(x_1,x_2)= ax_1^2 +\frac{r}{a}x_2^2+h_1^ix_1+h_2^ix_2+h_3^i, \qquad i=1,\dots,m,$$
such that
\begin{align}\label{gi-1h}
\begin{cases}q_i^r(\tilde{p}_{i-1,i}) &= q_{i-1}^s(\tilde{p}_{i-1,i})  \\
 \nabla q_i^r(\tilde{p}_{i-1,i}) &= \nabla q_{i-1}^s(\tilde{p}_{i-1,i})  
 \end{cases}
\end{align}
for  some $\tilde{p}_{i-1,i} \in (p_{i-1},p_i)$ and  with $(q_{i-1}^s+ \I_R)(x_1,x_2)\leq (q_i^r+\I_{R})(x_1,x_2)$ for all  $(x_1,x_2)\in \mathbb{R}^2$, and 
\begin{align}\label{gih}
   \begin{cases} q_i^r(\tilde{p}_{i,i}) &= q_{i}^s(\tilde{p}_{i,i}) \\
  \nabla q_i^r (\tilde{p}_{i,i}) &= \nabla q_{i}^s(\tilde{p}_{i,i})
  \end{cases}
\end{align}
for  some $\tilde{p}_{i,i} \in (\tilde{p}_{i-1,i},p_i)$ and with $(q_{i}^s+ \I_R)(x_1,x_2)\leq (q_i^r+\I_{R})(x_1,x_2)$ for all  $(x_1,x_2)\in \mathbb{R}^2$.

For simplicity, write $\tilde{p}_{i-1,i}= (t_1, y_{i-1,i})$ and $\tilde{p}_{i,i}= (t_1, y_{i,i})$. From the second equation in \eqref{gi-1h}, we have
$$\left(2at_2+h_1^i, 2\frac{r}{a}y_{i-1,i}+h_2^i\right)=\left(2at_2, 2\frac{s}{a}y_{i-1,i}+2\left(1-\frac{s}{a}\right)\left(-t_2+\left(\frac{2t_2}{m}\right)(i-1)\right)\right),$$
which implies $h_1^i=0$ and 
\begin{align*}
    y_{i-1,i}= \dfrac{h_2^i a}{2\left(s-r\right)}+\left(\dfrac{s-a}{s-r}\right)\left(-t_2+\left(\frac{2t_2}{m}\right)(i-1)\right).
\end{align*}

From the second equation in \eqref{gih}, it follows that 
$$\left(2at_2+h_1^i, 2\frac{r}{a}y_{i,i}+h_2^i\right)=\left(2at_2, 2\frac{s}{a}y_{i,i}+2\left(1-\frac{s}{a}\right)\left(-t_2+\left(\frac{2t_2}{m}\right)i\right)\right),$$
which implies
\begin{align*}
    y_{i,i}= \dfrac{h_2^i a}{2\left(s-r\right)}+\left(\dfrac{s-a}{s-r}\right)\left(-t_2+\left(\frac{2t_2}{m}\right)i\right)= y_{i-1}+\left(\dfrac{s-a}{s-r}\right)\frac{2t_2}{m}.
\end{align*}
Thus,
\begin{align}\label{yiyi-1}
    y_{i,i}-y_{i-1,i}=\left(\dfrac{s-a}{s-r}\right)\frac{2t_2}{m}
\end{align}
for every $i=1,\dots, m$, i.e., $\|\tilde{p}_{i-1,i}-\tilde{p}_{i,i}\|$ is a constant that does not depend on $i$.

We now compute the distance between $\tilde{p}_{i,i}$ and $\tilde{p}_{i,i+1}$. From the first equations in \eqref{gi-1h} and \eqref{gih}, we obtain
\begin{multline*}
    \dfrac{(r-s)}{a}y_{i-1,i}^2 +\left(h_2^i-2\left(1-\frac{s}{a}\right)\left(-t_2+\left(\frac{2t_2}{m}\right)(i-1)\right)\right)y_{i-1,i}\\
    +\left(1-\frac{s}{a}\right)\left(-t_2+\left(\frac{2t_2}{m}\right)(i-1)\right)^2+h_3^i=0,
\end{multline*}
and
\begin{align*}
    \dfrac{(r-s)}{a}y_{i,i}^2 +\left(h_2^i-2\left(1-\frac{s}{a}\right)\left(-t_2+\left(\frac{2t_2}{m}\right)i\right)\right)y_{i,i} +\left(1-\frac{s}{a}\right)\left(-t_2+\left(\frac{2t_2}{m}\right)i\right)^2+h_3^i=0.
\end{align*}
A straightforward computation yields
\begin{align*}
    h_2^i\left(\dfrac{s-a}{s-r}\right)=2\dfrac{(s-a)}{a}\left(-t_2+\left(\dfrac{2t_2}{m}\right)i\right)\left(1-\left(\dfrac{s-a}{s-r}\right)\right)+c,
\end{align*}
for every $i=1,\dots,m$, where $c$ is a constant that does not depend on $i$. Consequently,
$$h_2^{i+1}-h_2^{i}= 2\dfrac{(s-a)}{a} \left(\dfrac{2t_2}{m}\right)\left(\left(\dfrac{s-r}{s-a}\right)-1\right),$$
and
\begin{align}\label{yiyi+1}
y_{i,i+1}-    y_{i,i}=a\dfrac{h_2^{i+1}-h_2^{i}}{2(s-r)}=\dfrac{2t_2}{m}\left(1-\left(\dfrac{s-a}{s-r}\right)\right).
\end{align}

Consider the following sequence of functions
\begin{multline*}
u_m= (q_0^s+\I_{[-t_1,t_1]\times [-t_2, y_{0,1}]})\wedge (q_1^r+\I_{[-t_1,t_1]\times [ y_{0,1}, y_{1,1}]})\wedge (q_1^s+\I_{[-t_1,t_1]\times [ y_{1,1}, y_{1,2}]})\wedge\cdots\wedge\\
\wedge (q_m^r+\I_{[-t_1,t_1]\times [ y_{m-1,m}, y_{m,m}]})\wedge (q_m^s+\I_{[-t_1,t_1]\times [ y_{m,m}, t_2]}).
\end{multline*}
Since $y_{i,i}-y_{i-i,i}$ and $y_{i,i}-y_{i,i-1}$ do not depend on $i$, and $Z$ is  simple  and dually epi-translation invariant, we have
$$Z(u_m)=m Z(q_1^r+\I_{[-t_1,t_1]\times [ y_{0,1}, y_{1,1}]}) + mZ(q_1^s+\I_{[-t_1,t_1]\times [ y_{1,1}, y_{1,2}]})$$
and  this, combined with \eqref{mudanca_dominio}, implies
\begin{align}
\nonumber  Z(u_m)  &= \dfrac{m2t_1(y_{1,1}-y_{0,1})}{2t_12t_2}Z(q_1^r+\I_R)   + \dfrac{m2t_1(y_{1,2}-y_{1,1})}{2t_12t_2}Z(q_1^s+\I_R)\\
 \nonumber & =\dfrac{m2t_1}{2t_12t_2}\left(\dfrac{s-a}{s-r}\right)\left(\dfrac{2t_2}{m}\right)Z(q_1^r+\I_R) + \dfrac{m2t_1}{2t_12t_2}\left(1-\left(\dfrac{s-a}{s-r}\right)\right)\left(\dfrac{2t_2}{m}\right)Z(q_1^s+\I_R)\\
&  =\left(\dfrac{s-a}{s-r}\right)Z(q_1^r+\I_R)+ \left(1-\left(\dfrac{s-a}{s-r}\right)\right)Z(q_1^s+\I_R). \label{Z(um)}
\end{align}

Note that $Z(u_m)$ does not depend on $m$. To complete the proof that $\zeta$ is concave in the planar case, we use the facts that $u_m$ is $\tau$-convergent to $q^a+\I_R$ as $m\rightarrow \infty$, that $Z$ is $\tau$-upper semicontinuous, and that \eqref{Z(um)}  and \eqref{mudanca_zeta} hold,  so that
\begin{align*}
    \zeta(2^2a)V_2(R)&=Z(q^a+\I_R)\\
    & \geq \limsup_{m\rightarrow \infty} Z(u_m)\\
    &= \limsup_{m\rightarrow \infty} \left(\left(\dfrac{s-a}{s-r}\right)Z(q_1^r+\I_R)+ \left(1-\left(\dfrac{s-a}{s-r}\right)\right)Z(q_1^s+\I_R)\right)\\
    &=V_2(R)\left(\left(\dfrac{s-a}{s-r}\right) \zeta(2^2r)+ \left(1-\left(\dfrac{s-a}{s-r}\right)\right)\zeta(2^2s)\right).
\end{align*}
Setting $\lambda= \dfrac{s-a}{s-r}$, since $0\leq s<a<r$, we have  $0<\lambda<1$ and 
$$\zeta(2^2\lambda r+2^2(1-\lambda)s)\geq \lambda \,\zeta(2^2r)+(1-\lambda)\,\zeta(2^2s).$$
Since $r$ and $s$ were arbitrarily chosen, we arrive at
$$\zeta(\lambda r+(1-\lambda)s)\geq \lambda\, \zeta(r)+(1-\lambda)\,\zeta(s)$$
for every $\lambda\in (0,1)$, i.e., $\zeta$ is concave.
\medskip

To prove that $\zeta$ is concave for $n\geq 3$, we extend the quadratic functions $q^a, q_i^s$ and $q_i^r$ so as to obtain quadratic functions on $\mathbb{R}^n$: the planar graphs remain unchanged, and in  the directions of the $x_j$-axes, $j=3, \dots, n$, we extend the domain by taking the Cartesian product with the intervals $[-t_j,t_j]$, so that the resulting domain is
$
R = [-t_1,t_1]\times \cdots \times [-t_n,t_n].
$
Thus, 
\begin{align*}
    q^a(x)= \langle  x, \diag(a,1, \dots, 1)x\rangle,
\end{align*}
the quadratic function $q_i^s$ is given by 
\begin{multline*}
    q_i^s(x)= \langle x, \diag (a, a^{-1}s, 1,\dots, 1)x\rangle + 2\left(1-\frac{s}{a}\right)\left(-t_2+\left(\frac{2t_2}{m}\right)i\right)x_2\\
   \nonumber -\left(1-\frac{s}{a}\right)\left(-t_2+\left(\frac{2t_2}{m}\right)i\right)^2,
\end{multline*}
where $x= (x_1,x_2, \dots, x_n)$, and the quadratic function $q_i^r$ is obtained from 
\begin{align*}
    \tilde{q^r_i}(x)= \langle x, \diag(a, a^{-1}r, 1,\dots, 1)x\rangle.
\end{align*}
by a rigid rotation in the $(x_1,x_2)$-plane that preserves the coordinate axes and fixes the remaining coordinates $x_3,\dots,x_n$.

Now,
$$p_i=\left(t_1, -t_2+\left(\frac{2t_2}{m}\right)i, 1,\dots, 1\right),$$
for every $i=0, \dots, m$, and  $R= [-t_1,t_1]\times \cdots \times [-t_n,t_n]$. We have 
$$q_i^s(p_i)=q^a(p_i),\, \nabla q_i^s(p_i)= \nabla q^a(p_i) \ \mbox{ and } (q_i^s+\I_R)(x)  \leq (q^a+\I_R)(x) \ \mbox{for all } x\in \mathbb{R}^n.$$
The epi-graphs of $q_i^r+\I_R$ are contained in $\epi(q_{i-1}^s+\I_R)$ and in  $\epi(q_{i}^s+\I_R)$. They coincide with $q_{i-1}^s$ along the hyperplane $x_2=0$, translated by $\tilde{p}_{i-1,i}$, and  with $q_{i}^s$ along the hyperplane $x_2=0$, translated by $\tilde{p}_{i,i}$. As before, write 
$$\tilde{p}_{i-1,i}= (t_1, y_{i-1,1}, t_3, \dots, t_n),\, \tilde{p}_{i,i}= (t_1, y_{i,i}, t_3,\dots, t_n),$$ 
and define 
\begin{multline*}
u_m= (q_0^s+\I_{[-t_1,t_1]\times [-t_2, y_{0,1}]\times J })\wedge (q_1^r+\I_{[-t_1,t_1]\times [ y_{0,1}, y_{1,1}]\times J})\wedge (q_1^s+\I_{[-t_1,t_1]\times [ y_{1,1}, y_{1,2}]\times J})\wedge\cdots\wedge\\
\wedge (q_m^r+\I_{[-t_1,t_1]\times [ y_{m-1,m}, y_{m,m}]\times J})\wedge (q_m^s+\I_{[-t_1,t_1]\times [ y_{m,m}, t_2]\times J}),
\end{multline*}
where $J=[-t_3,t_3]\times \cdots \times [-t_n,t_n]$. 
\medskip

As $m\rightarrow \infty$, the sequence $u_m$ is $\tau$-convergent to $q^a+\I_R$. For the  general case, we also have  equations \eqref{yiyi-1} and \eqref{yiyi+1}. Therefore, analogously to the two-dimensional case,
\begin{align*}
    \zeta(2^na)V_n(R)&=Z(q^a+\I_R)\\
    & \geq \limsup_{m\rightarrow \infty} Z(u_m)\\
    &=V_n(R)\left(\left(\dfrac{s-a}{s-r}\right) \zeta(2^nr)+ \left(1-\left(\dfrac{s-a}{s-r}\right)\right)\zeta(2^ns)\right),
\end{align*}
and hence $\zeta$ is a concave function.
\end{proof}

\begin{lema}\label{concave1}
The function $\zeta$ given by \eqref{zeta} belongs to $\conc$.
\end{lema}

\begin{proof}
By Lemma \ref{conccave}, $\zeta$ is concave. Note that $a^{\frac{1}{n}}q+\I_{C}\in \LP\RR$ is $\tau$-convergent to  $\I_C$ as $a\rightarrow 0$. By Lemma \ref{positive}, and using that $Z$ vanishes on  indicator functions, we get
$$0\leq \limsup_{a\rightarrow 0}\zeta(a) = \limsup_{a\rightarrow 0} \dfrac{Z(a^{\frac{1}{n}}q+\I_C)}{V_n(C)}\leq \dfrac{Z(\I_C)}{V_n(C)}=0.$$

Consider now the sequence of functions
$$u_k(x)= \left\langle x, \diag(k2^{-n},1, \dots, 1)x\right\rangle $$
for $k\in\mathbb{N}$.  Note that $u_k+\I_{[0,1/k]\times[0,1]^{n-1}} $ is $\tau$-convergent to $\left\langle \cdot, \diag(0,1, \dots, 1)\cdot\right\rangle+\I_{\{0\}\times[0,1]^{n-1}}$ as $k\rightarrow \infty$. Since $Z$ is a $\tau$-upper semicontinuous and simple valuation, we obtain
$$\limsup_{k\rightarrow \infty} Z(u_k+\I_{[0,1/k]\times[0,1]^{n-1}})=0.$$
On the other hand, by formulas  \eqref{mudanca_dominio} and \eqref{mudanca_zeta}, we have
\begin{align*}
    Z(u_k+\I_{[0,1/k]\times[0,1]^{n-1}})&= Z\left(\langle \cdot, \diag(k2^{-n},1, \dots, 1)\cdot\rangle
    +\I_{[0,1/k]\times [0,1]^{n-1}}\right)\\
    &= \dfrac{Z\left(\langle \cdot, \diag(k2^{-n},1, \dots, 1)\cdot\rangle
    +\I_{C}\right)}{V_n(C)}\dfrac{1}{k}\\
    &= \dfrac{\zeta(k)}{k}.
\end{align*}
Hence
$$\limsup_{k\rightarrow \infty} \dfrac{\zeta(k)}{k}=0.$$
This concludes the proof.
\end{proof}

Next, we show that the function $\zeta$ completely determines the valuation $Z$. To this end, let $P\LQ\RR$ denote the family of convex functions of the form
$$u=\bigwedge_{i=1}^m w_i,$$
where  each $w_i$ is either a piecewise linear-quadratic function with polytopal domain or a cylinder function with  polytopal domain. Let $Z$ be a simple, $\tau$-upper semicontinuous, equi-affine and dually epi-translation invariant valuation on $\LP\RR$ that vanishes on indicator functions of polytopes and cylinder functions, and let $\zeta$ be defined from $Z$ by \eqref{mudanca_zeta}. Then the value of $Z$ is well defined on the constituent functions appearing in the representation of $u \in P\LQ\RR$, namely piecewise linear-quadratic functions with polytopal domain and cylinder functions with polytopal domain. In particular, its values on quadratic components are determined by $\zeta$, and hence $Z$ is determined on all of $P\LQ\RR$. Note that  $\PA\RR\cap \LP(\mathbb{R}^n)\subset P\LQ\RR$ and, since $\PA\RR\cap \LP(\mathbb{R}^n)$ is dense in $\LP\RR$ (see Lemma \ref{piecewise}),  every $u\in \LP\RR$ can be approximated by elements of   $P\LQ\RR$. The $\tau$-upper semicontinuity of $Z$ then implies that for any sequence $u_k\in P\LQ\RR$, $k\in\mathbb{N}$, with $u_k\seq u$, we have 
\begin{align}\label{1234}
  Z(u)\geq \limsup_{k\rightarrow \infty} Z(u_k).  
\end{align}
To show that $Z$ is uniquely determined by $\zeta$, it suffices to construct, for every  $u\in \LP(\mathbb{R}^n)$, a sequence  $u_k\in P\LQ(\mathbb{R}^n)$, $k\in\mathbb{N}$, such that equality holds in \eqref{1234}, i.e.,
\begin{align}\label{equality}
    Z(u)=\sup\Bigl\{\limsup_{k\rightarrow \infty} Z(u_k)\mid u_k\in P\LQ(\mathbb{R}^n),\, u_k\seq u\Bigr\}.
\end{align}
Since each $Z(u_k)$ is determined by $\zeta$, it follows that the right-hand side of \eqref{equality} is completely determined by $\zeta$.

We will show that it is enough to prove \eqref{equality} for the special class of functions $u_{\lambda,\mu}$ introduced in Subsection~\ref{bauschke2011convex}.  
For this, we require the following lemma.

\begin{lema}\label{invariant}
Let $Z:\LP\RR\rightarrow\mathbb{R}$ be a simple, $\tau$-upper semicontinuous, equi-affine and  dually epi-translation   invariant valuation that   vanishes on  indicator functions of polytopes and cylinder functions. Let $u \in \LP(\mathbb{R}^n)$ and let $v$ be an affine function on $\mathbb{R}^n$.
Then we have
$$
Z(u\,\square\, (v+\I_J))=Z(u),
$$
where $J\subset \mathbb{R}^n$ is a closed line segment. 
\end{lema}

\begin{proof}
We begin by considering line segments $J\subset \mathbb{R}^n$ in a fixed line $L$. 
Since $Z$ is a translation invariant valuation, $Z(u\,\square\, (v+\I_J))$ depends only on the length $s=V_1(J)$. Thus, the function $g_u: [0,\infty)\rightarrow [0, \infty)$
defined by
\begin{align*}
  g_u(s)= Z(u\,\square\, (v+\I_J)), \quad  s=V_1(J),
\end{align*}
is well-defined. As in the proof of Lemma \ref{cylinder}, let $J_1$ and $J_2$ be line segments lying on the same line, with disjoint relative interiors, and suppose that their intersection consists of a single point $p$ such that $v(p)=0$. Then
\begin{align*}
    (u\,\square\, (v+\I_{J_1}))\wedge (u\,\square\, (v+\I_{J_2})) & =u\,\square\, (v+\I_{J_1\cup J_2})\\
    (u\,\square\, (v+\I_{J_1}))\vee (u\,\square\, (v+\I_{J_2})) & =u.
\end{align*}

Let $s_1=V_1(J_1)$ and $s_2=V_1(J_2)$. Since $Z$ is  a valuation, it follows that
\begin{align*}
    Z(u\,\square\, (v+\I_{J_1\cup J_2}))+Z(u)= Z(u\,\square\, (v+\I_{J_1}))+Z(u\,\square\, (v+\I_{J_2})),
\end{align*}
which implies
\begin{align*}
    g_u(s_1+s_2)+g_u(0)=g_u(s_1)+g_u(s_2)
\end{align*}
for every $s_1,s_2\geq 0$. Hence, $g_u(s)-g_u(0)$ satisfies Cauchy’s functional equation and depends upper semicontinuously on $s$. Therefore, there exists a real-valued function $u \mapsto Y(u)$ such that
\begin{align*}
    g_u(s)=Y(u)\, s+g_u(0),
\end{align*}
for every $s\geq 0$. In particular,
\begin{align}\label{1.18}
Z(u\,\square\, (v+\I_J))  = Y(u)\,V_1(J)+Z(u)  
\end{align}
holds for every $u\in\LP\RR$, every  affine function $v$ on $\mathbb{R}^n$, and every closed line segment $J\subset L$.

For $J$ sufficiently long in $L$, and by definition of the epi-sum, we can find closed half-spaces $H_1^+$ and $H_2^-$,  orthogonal to $J$,  such that   
$$u\,\square\, (v+\I_{J})+\I_{(\dom u+ J)\cap H_1^+\cap H_2^-}$$ 
is a non-degenerate cylinder function $w\in \LP\RR$. Set $J_1=J\cap H_1^+$ and $J_2=J\cap H_2^-$. Note that 
\begin{multline*}
    Z(u\,\square\, (v+\I_{J}))\\
    = Z(u\,\square\, (v+\I_{J})+\I_{\dom (u\,\square\, (v+\I_{J}))\cap H_1^+})+ Z(w)+Z(u\,\square\, (v+\I_{J})+\I_{\dom (u\,\square\, (v+\I_{J}))\cap H_2^-}).
\end{multline*}
Since $Z$ is a simple dually epi-translation invariant valuation that vanishes on cylinder functions,
\begin{align*}
    Z(u\,\square\, (v+\I_{J}))= Z(u\,\square\, (v+\I_{I_1+I_2})).
\end{align*}
Using this together with \eqref{1.18}, we obtain
$$ Y(u)\,V_1(J)= Y(u)\,V_1(J_1+J_2)= Y(u)\,(V_1(J_1)+V_1(J_2)),$$
and, since $V_1(J)>V_1(J_1)+V_1(J_2)$, it follows that
$$Y(u)=0.$$
Substituting into \eqref{1.18}, we conclude that 
$$Z(u\,\square\, (v+\I_{J}))= Z(u)$$
for every $u\in \LP\RR$,  and thus the proof is complete.
\end{proof}

Let $\lambda,\mu>0$ and $u\in \LP(\mathbb{R}^n)$. Recall from \eqref{ulm} that
\begin{equation}
     u_{\lambda,\mu}= u\,\square\, \Big(\frac{\lambda}{2}\|\cdot\|^2+\I_{\mu C}\Big).
\end{equation}
Let $\tilde{w}_k\in \LP(\mathbb{R}^n)$, $k \in \mathbb{N}$, be a sequence of functions given by Lemma \ref{1-dim}, i.e.,
$$\tilde{w}_k=  w_k^{(1)} \, \square \cdots \square \,  w_k^{(n)} = \lfloor J_1^{k,1} \rfloor \,\square \cdots \square\, \lfloor J_{m_k^{(1)}}^{k,1} \rfloor \,\square \cdots \square\, \lfloor J_1^{k,n} \rfloor \,\square \cdots \square\, \lfloor J_{m_k^{(n)}}^{k,n} \rfloor,$$
and $u\,\square\, \lambda \tilde{w}_k \seq u_{\lambda,\mu}$ as  $k\to \infty$.

Then, by Lemma \ref{invariant}, if $Z:\LP\RR\rightarrow\mathbb{R}$ is a simple, $\tau$-upper semicontinuous, equi-affine and  dually epi-translation   invariant valuation that   vanishes on  indicator functions of polytopes and cylinder functions, we obtain
\begin{align}\label{bigger}
Z\bigl(u_{\lambda,\mu}\bigr)\geq \limsup_{k\rightarrow \infty} Z\left(u\,\square\, \lambda \tilde{w}_k\right)=  Z(u). 
\end{align}

Now suppose, by contradiction, that \eqref{equality} does not hold for some $u\in \LP\RR$, i.e.,
$$
Z(u)>\sup\bigl\{\limsup_{k\rightarrow \infty} Z(v_k): v_k\in P\LQ\RR,\ v_k\seq u\bigl\}.
$$
Then, by Theorem \ref{equivalence}, there exist $\rho>0$  and $\tilde{\mu} >0$ such that
\begin{align}\label{equa21}
    Z(u)> \limsup_{k\rightarrow \infty} Z(v_k)+\rho \, V_n(\dom u)
\end{align}
for every sequence $v_k\in P\LQ\RR$, $k\in\mathbb{N}$, whose Lipschitz constants in 
$\interior(\dom v_k)$ are uniformly bounded by some $L > 0$, and such that 
$v_k$ converges uniformly to $u$ on every compact set $C$ that does not contain 
boundary points of $\dom u$ and also does not contain boundary points of $\dom u_{\lambda,\mu}$ for $\mu \le \tilde{\mu}$.

Let $v_{k,\mu}\in P\LQ\RR$ be such that $v_{k,\mu}\seq u_{\lambda,\mu}$ as $k\to \infty$. Then the Lipschitz constants  of $v_{k,\mu}$ are uniformly bounded and  $v_{k,\mu}$ epi-converges to $u_{\lambda,\mu}$ as $k\to \infty$. By Theorem \ref{equivalence}, $v_{k,\mu}$ converges uniformly to $u_{\lambda,\mu}$ on every compact set that does not contain boundary points of $\dom u_{\lambda,\mu}$.  In particular, let $D$ be a compact set  such that does not contain boundary points 
of $\dom u_{\lambda,\mu}$ and also does not contain boundary points of $\dom u$.  Since $u_{\lambda,\mu}$ epi-converges to $u$ as $\mu \to 0$ and $D$ does not 
contain boundary points of $\dom u$, Theorem~\ref{equivalence} implies that 
for every $\varepsilon > 0$ there exists $\delta > 0$ such that
\begin{align*}
    |u(x)-u_{\lambda,\mu}(x)|< \frac{\varepsilon}{2}
\end{align*}
for all $\mu\le \delta$ and $x\in D$. Similarly, since $v_{k,\mu}$ converges uniformly to $u_{\lambda,\mu}$ on $D$,  there exists $k_0 > 0$ such that
\begin{align*}
    |u_{\lambda,\mu}(x)- v_{k,\mu}(x)|< \frac{\varepsilon}{2}
\end{align*}
for all $k\geq k_0$ and $x\in D$. Taking $\mu \le \min\{\delta,\tilde{\mu}\}$ and $k\geq k_0$, we obtain 
\begin{align*}
    |u(x)-v_{k,\mu}(x)|\leq |u(x)-u_{\lambda,\mu}(x)|+ |u_{\lambda,\mu}(x)- v_{k,\mu}(x)|< \varepsilon
\end{align*}
for all $x\in D$. Thus, $v_{k,\mu }$ converges uniformly to $u$ on every compact set  $D$ that does not contain boundary points of $\dom u_{\lambda,\mu} $ and also does not contain boundary points of $\dom u$ for $\mu\leq \tilde{\mu}$. Then, by  \eqref{bigger} and  \eqref{equa21}, we have 
\begin{align*}
 Z\left(u_{\lambda,\mu}\right)>  \limsup_{k\rightarrow \infty} Z(v_{k,\mu})+\rho \, V_n(\dom u)   
\end{align*}
for every $v_{k,\mu}\in P\LQ\RR$ with $v_{k,\mu}\seq u_{\lambda,\mu}$ as $k\to \infty$, provided that $\mu$ is sufficiently small.   Since $V_n$ depends continuously on the domain and $\dom u_{\lambda,\mu}= \dom u +\mu C$, it follows that
\begin{align*}
Z\left(u_{\lambda,\mu}\right)>  \limsup_{k\rightarrow \infty} Z(v_{k,\mu})+\frac{\rho}{2}\, V_n(\dom u_{\lambda,\mu})    
\end{align*}
for  sufficiently small $\mu$ and every $v_{k,\mu}\in P\LQ\RR$ with   $v_{k,\mu}\seq u_{\lambda,\mu}$ as $k\to \infty$. This shows that if \eqref{equality} fails for some $u\in \LP\RR$, then it  also fails for  $u_{\lambda,\mu}$ for sufficiently small $\mu>0$.

It therefore remains to establish \eqref{equality} for the functions $u_{\lambda,\mu}$. The proof of this fact requires several auxiliary results, which we state below and prove in Appendix~\ref{apendice}.

Let $l^i$ and $l^c$ be piecewise affine functions in $\LP\RR$ such that 
\begin{align}\label{li_lc}
l^c(x)< u(x)< l^i(x)
\end{align}
for all $x\in\mathbb{R}^n$. We will show that for every choice of such $u$, $l^i$, $l^c$, and every $\rho>0$, one can construct a function $v\in P\LQ(\mathbb{R}^n)$ satisfying
$$
\dom l^i \subset \dom v \subset \dom l^c
$$
and
\begin{align}\label{eq10.22}
    Z(u)\leq Z(v)+\rho \, V_n(\dom u) \quad \mbox{and} \quad  l^c(x)\leq v(x)\leq l^i(x) \text{ for all } x\in\mathbb{R}^n.
\end{align}
This shows  that for any given $u$, there exists 
$v\in P\LQ\RR$ arbitrarily close to  $u$ such that $Z(v)$ is nearly as large as $Z(u)$, thereby proving \eqref{equality}. From now on, we assume that \eqref{li_lc} holds and that $\rho>0$.

We note that the arguments leading to the construction of $v$ are somewhat involved and require several technical steps. For readability, the auxiliary results needed in this part of the proof are stated in the main text below and proved in Appendix~\ref{apendice}.

Let $N \subset \dom u$ denote the set of interior points where $u$ is twice differentiable. 
By Aleksandrov's theorem \cite{aleksandrov1939second}, we have $V_n(N) = V_n(\dom u)$.  We now proceed to establish the preliminary results needed for the construction.

\begin{prop}\label{regular}
Let $Z:\LP\RR\rightarrow\mathbb{R}$ be a simple, $\tau$-upper semicontinuous, equi-affine and  dually epi-translation   invariant valuation that   vanishes on  indicator functions of polytopes and cylinder functions. Let $u \in \LP(\mathbb{R}^n)$,
let $x_0\in N$ be such that $\det\D u(x_0)>0$,  let $\rho>0$, and let $l^i, l^c$  be the functions satisfying \eqref{li_lc}. Then there exist $r(x_0)>0$, 
polytopes $P_r(x_0)$ containing $x_0$ in their interiors, and  functions $v_r^{x_0}\in P\LQ\RR$ such that
\begin{align}\label{eq10.47}
  Z(u+\I_{P_r(x_0)}) \leq Z(v_r^{x_0}+\I_{P_r(x_0)}) + \frac{\rho}{2}V_n(P_r(x_0))
\end{align}
and
\begin{align}\label{eq10.48}
l^c(x)\leq v_r^{x_0}(x)\leq l^i(x)    
\end{align}
for all $x\in \mathbb{R}^n$ and $0<r\leq r(x_0)$. Moreover, for $0<r\leq r(x_0)$,  the function $v_r^{x_0}$ can be chosen so such that
$u$ and $v_r^{x_0}$ share a common supporting hyperplane at $x_0$, and 
\begin{align}\label{eq4111}
 \{x\in \mathbb{R}^n\mid v_r^{x_0}(x)\geq u(x)\}   \subset P_r(x_0), 
\end{align}
while the  polytopes $P_r(x_0)$ satisfy 
\begin{align}\label{eq10.26..}
\diam(P_r(x_0))^n= n^{\frac{n}{2}}\, V_n(P_r(x_0)).    
\end{align}
\end{prop}

To prove Proposition \ref{regular}, we compare $u$ in a neighborhood of $x_0$ relative to $\dom u$, with  a suitable restriction of 
the quadratic function $v_r^{x_0}+\I_{P_r(x_0)}$ that shares a common tangent hyperplane with $u$ at $x_0$. We then show that $Z(v_r^{x_0}+\I_{P_r(x_0)})$  is nearly as large as $Z(u+\I_{P_r(x_0)})$. However, if $x_0\in N$ is a point where $\det \D u(x_0)=0$, we instead compare $u$ with piecewise affine functions.

\begin{prop}\label{regular2}
Let $Z:\LP\RR\rightarrow\mathbb{R}$ be a simple, $\tau$-upper semicontinuous, equi-affine and  dually epi-translation   invariant valuation that   vanishes on  indicator functions of polytopes and cylinder functions. Let $u \in \LP(\mathbb{R}^n)$,
let $x_0\in N$ be such that $\det \D u(x_0)=0$,  let $\rho>0$, and let $l^i, l^c$  be the functions satisfying \eqref{li_lc}. Then there exist $r(x_0)>0$ and   polytopes $P_r(x_0)$ containing $x_0$ in their interiors  such that
\begin{align}\label{eq10.67}
  Z(u+\I_{P_r(x_0)})\leq \frac{\rho}{2}\, V_n(P_r(x_0))
\end{align}
and
\begin{align}\label{eq10.66}
 l^c(x)\leq v_r^{x_0}(x)\leq l^i(x)    
\end{align}
for all  $x\in \mathbb{R}^n$ and $0<r\leq r(x_0)$, where
\begin{align}\label{eq41112}
  v_r^{x_0}= l_{x_0}\vee l^c,  
\end{align}
and $l_{x_0}$ is the affine function defining the tangent hyperplane to $u$ at $x_0$, coinciding with $u$ at $x_0$. Moreover, for $0<r\leq r(x_0)$, the 
polytopes $P_r(x_0)$ satisfy 
\begin{align}\label{eq10.53..}
\diam(P_r(x_0))^n= n^{\frac{n}{2}}\, V_n(P_r(x_0)).    
\end{align}
\end{prop}

The following proposition establishes an ``absolute continuity'' property of $Z$ for functions of the form
\[
u=\tilde{u}\,\square\,\frac{\lambda}{2}(\|\cdot\|^2+\I_{\mu C}),
\]
where $\tilde{u}\in\LP(\mathbb{R}^n)$ and $\lambda,\mu>0$.

\begin{prop}\label{singular}
Let $Z:\LP\RR\rightarrow\mathbb{R}$ be a simple, $\tau$-upper semicontinuous, equi-affine and  dually epi-translation   invariant valuation that   vanishes on  indicator functions of polytopes and cylinder functions. Let $u = \tilde{u}\,\square\, \frac{\lambda}{2}\big(\|\cdot\|^2+\I_{\mu C}\big)$, where $\tilde{u}\in \LP(\mathbb{R}^n)$,  $ \lambda>0$, and $0<\mu<1$.   There is a constant $c_{\lambda}$, depending on the dimension $n$ and on $\lambda$, such that
\begin{align*}
    Z(u+\I_{\dom u \, \cap \,P})\leq c_{\lambda}\, V_n(\dom u\cap P) 
\end{align*}
for every polytope $P\subset \mathbb{R}^n$.
\end{prop}

Combining Proposition~\ref{regular}, Proposition~\ref{regular2}, and Proposition~\ref{singular}, we are now able to establish \eqref{equality} for the functions $u_{\lambda,\mu}$. By the reduction established earlier in this section, this implies \eqref{equality} for every function in $\LP\RR$.

\begin{prop}\label{prop10.3}
Let $Z:\LP\RR\rightarrow\mathbb{R}$ be a simple, $\tau$-upper semicontinuous, equi-affine and  dually epi-translation   invariant valuation that   vanishes on  indicator functions of polytopes and cylinder functions. For every $u = \tilde{u}\,\square\, \frac{\lambda}{2}\big(\|\cdot\|^2+\I_{\mu C}\big)\in \LP\RR$, where $\tilde{u}\in \LP(\mathbb{R}^n)$,  $ \lambda>0$, and $0<\mu<1$, we have 
   $$Z(u)=\sup\Bigl\{\limsup_{k\rightarrow\infty} Z(v_k)\mid v_k\in P\LQ\RR, v_k\seq u   \Bigr\}.$$ 
\end{prop}

\begin{proof}
By Proposition~\ref{regular}, Proposition~\ref{regular2}, and Proposition~\ref{singular}, we construct a function
$v\in P\LQ\RR$ satisfying \eqref{eq10.22}. As before, let $N$ be the set of  points in $\interior(\dom u)$ where $u$ is twice differentiable,  and let $\mathcal{L}$ be the collection of sets  $P_r(x)$
for $x\in N$ and $0<r\leq r(x)$, where $P_r(x)$ and $r(x)$ are defined in  Proposition \ref{regular} and Proposition \ref{regular2}.  By \eqref{eq10.26..} and \eqref{eq10.53..}, the collection $\mathcal{L}$ forms a regular family for $N$, and this remains true if we restrict to those sets with  $V_n(P_r(x))\leq \delta$  for some $\delta>0$.

Now fix $\eta>0$ such that
\begin{align*}
\eta \leq \delta \qquad \mbox{and} \qquad     \eta\leq \dfrac{\rho}{2c_{\lambda}} \,V_n(\dom u),
\end{align*}
where $c_{\lambda}$ is given by Proposition \ref{singular}. 
By Vitali's Theorem \ref{VCT}, there exist pairwise disjoint sets $P_{r_1}(x_1), \dots, P_{r_m}(x_m)\in\mathcal{L}$ such that
\begin{align*}
    V_n(\dom u)-\sum_{i=1}^m V_n(P_{r_i}(x_i))\leq \eta.  
\end{align*}

Let  $v_{r_1}^{x_1}, \dots, v_{r_m}^{x_m}$ be the functions from  Proposition \ref{regular} and Proposition \ref{regular2}.  By \eqref{eq4111} and \eqref{eq41112}, it follows that for $i\neq j$, the functions $v_{r_i}^{x_i}$ and $v_{r_j}^{x_j}$ do not coincide at points where
$v_{r_i}^{x_i}> u$ and $v_{r_j}^{x_j}> u$.
Thus, we can choose a piecewise affine function $l$ on $\mathbb{R}^n$ such that
$$ \qquad l(x)\leq u(x) \quad \mbox{ for all }  x\in\dom u,$$
and such that, for every pair $i \neq j$, there exist polytopes $\tilde{P}_{i} \subset \tilde{P}_{ij} \subset \dom u$ with 
$$
v_{r_i}^{x_i} = l + \I_{\dom u} \text{ on } \bd \tilde{P}_{i}, \quad 
v_{r_i}^{x_i} = v_{r_j}^{x_j} \text{ on } \bd \tilde{P}_{ij}.
$$
Here, $\tilde{P}_i \subset \tilde{P}_{ij}$ ensures that $v_{r_i}^{x_i}$ coincides with $l + \I_{\dom u}$  before it coincides with $v_{r_j}^{x_j}$, which allows us to control the order of interactions and avoid undesired intersections between these functions.

We then define 
$$v= \bigvee_{i=1}^m v_{r_i}^{x_i}\vee (l+\I_{\dom u}).$$
This construction ensures that $v\in P\LQ\RR$, and   \eqref{eq10.48} and \eqref{eq10.66} imply that $ l^c(x) \leq v(x) \leq l^i(x)$ for all $x \in \mathbb{R}^n$.

Next, we partition
$$\dom u\setminus \bigcup_{i=1}^mP_{r_i}(x_i)$$
into polytopes $Q_1, \dots, Q_k$ so that
$$\dom u= \bigcup_{i=1}^mP_{r_i}(x_i)\cup \bigcup_{s=1}^k Q_s.$$
Since the elements $P_{r_1}(x_1), \dots, P_{r_m}(x_m)$ of $\mathcal{L}$ are pairwise disjoint, the intersection of any two distinct sets from the 
collection $\{P_{r_i}(x_i)\}_i\cup~\{Q_s\}_s$   is either empty or an at most $(n-1)$-dimensional polytope contained in the domain of $ u$.   Since $Z$ is simple, we obtain
\begin{align*}
    Z(u)&= \sum_{i=1}^m Z(u+\I_{P_{r_i}(x_i)})+\sum_{s=1}^k Z(u+\I_{ \dom u\,\cap\, Q_s}).
\end{align*}

By construction, if $x_i$ is a point where $u$ is twice differentiable and  $\det \D u(x_i)=0$, then $v_{r_i}^{x_i}$ is a piecewise affine function in $P\LQ\RR$.  
If instead $\det\D u(x_i)>0$, then $v_{r_i}^{x_i}$ consists of a quadratic function (whose epi-graph is contained in $\epi (u)$), together with pieces of 
cylinder and piecewise affine functions.  
Since $Z$ vanishes on cylinder and piecewise affine functions,
$$
Z(v+\I_{P_{r_i}(x_i)})= Z(v_{r_i}^{x_i}+\I_{P_{r_i}(x_i)}).
$$
Furthermore, since each $v+\I_{Q_s}$ is either a cylinder function 
or a piecewise affine function for $s = 1, \dots, k$, it follows that
$$
Z(v)= \sum_{i=1}^m Z(v_{r_i}^{x_i}+\I_{P_{r_i}(x_i)}).
$$

Finally, applying Proposition \ref{regular},  Proposition \ref{regular2}, and Proposition \ref{singular}, we conclude that
\begin{align*}
    Z(u)&= \sum_{i=1}^m Z\bigl(u+\I_{P_{r_i}(x_i)}\bigr)+\sum_{s=1}^k Z(u+\I_{ \dom u\,\cap \,Q_s})\\
    &\leq \sum_{i=1}^m\left( Z\bigl(v_{r_i}^{x_i}+\I_{P_{r_i}(x_i)}\bigr)+\dfrac{\rho}{2}V_n(P_{r_i}(x_i))\right)+c_{\lambda}\sum_{s=1}^k V_n(\dom u\cap Q_s)\\
    &\leq  Z(v) + \sum_{i=1}^m\dfrac{\rho}{2}V_n(P_{r_i}(x_i)) +c_{\lambda}\eta\\
    &\leq  Z(v) +  \dfrac{\rho}{2}V_n(\dom u)+c_{\lambda}\dfrac{\rho}{2c_{\lambda}} V_n(\dom u)\\
    &=  Z(v) +\rho \, V_n(\dom u).
\end{align*}
Since $\rho>0$ is arbitrary,  this proves the proposition. 
\end{proof}

Proposition~\ref{prop10.3} establishes \eqref{equality}, thereby completing the argument that the function $\zeta$ uniquely determines the valuation $Z$. 
On the other hand, by \cite[Theorem~1]{baetaludwigsemicontinuity} and \cite[Lemma~12]{baeta1d}, for every $\zeta\in\conc$ the functional
\[
Z_\zeta(u)=\int_{\mathbb{R}^n}\zeta(\det \D u(x))\,\mathrm{d}x
\]
is a simple, $\tau$-upper semicontinuous, equi-affine and dually epi-translation invariant valuation on $\LP(\mathbb{R}^n)$ that vanishes on indicator functions of polytopes and cylinder functions, and satisfies
\[
Z_\zeta(aq+\I_P)=\zeta(a)V_n(P)
\]
for every $a\ge0$ and every polytope $P\subset\mathbb{R}^n$. We therefore obtain the following existence and uniqueness result.

\begin{prop}\label{uniquen}
For a given $\zeta \in \conc$, there exists a unique functional 
$$
Z: \LP(\mathbb{R}^n)\rightarrow\mathbb{R}
$$ 
satisfying the following properties:
\begin{itemize}
    \item[(i)] $Z$ is $\tau$-upper semicontinuous;
    \item[(ii)] $Z$ is a simple, equi-affine and dually epi-translation invariant valuation; 
    \item[(iii)] $Z$  vanishes on indicator functions of polytopes and cylinder functions;
    \item[(iv)] $Z(aq+\I_P)= \zeta(a)V_n(P)$ for every $a \geq 0$ 
    and  polytope $P$, where $q(x)=\frac{1}{2}\sum_{i=1}^n x_i^2$.
\end{itemize}
\end{prop}

\begin{proof}
Let $Z:\LP\RR\rightarrow\mathbb{R}$ be a valuation satisfying properties $(i)-(iv)$.  By property $(iv)$, together with the fact that $Z$ is simple, dually epi-translation invariant, and vanishes on indicator functions of polytopes and cylinder functions, the value of $Z$ on every function in $P\LQ\RR$ is completely determined by $\zeta$. Since every $u\in\LP\RR$ can be approximated in the $\tau$-sense by functions in $P\LQ\RR$, and $Z$ is $\tau$-upper semicontinuous, it follows from the reduction established before Proposition~\ref{prop10.3} and Proposition~\ref{prop10.3} that the value of $Z(u)$ is determined by its values on $P\LQ\RR$, and hence by $\zeta$. This proves that $Z$ is uniquely determined by $\zeta$.
\end{proof}

Finally, we are ready to complete the proof of the case where
$Z:\LP\RR\rightarrow\mathbb{R}$ is a simple, $\tau$-upper semicontinuous,
equi-affine and dually epi-translation invariant valuation that vanishes on
indicator functions of polytopes and cylinder functions.

\begin{prop}\label{prop}
Let $Z:\LP\RR\rightarrow\mathbb{R}$ be a simple, $\tau$-upper semicontinuous, equi-affine and  dually epi-translation   invariant valuation that   vanishes on  indicator functions of polytopes and cylinder functions. Then there is a   function $\zeta \in\conc$  such that
\begin{align*}
    Z(u)= \int_{\dom u}\zeta(\det \D u(x))\dif x
\end{align*}
for every $u\in \LP\RR$.
\end{prop}

\begin{proof}
Let $\zeta\in \conc$ be the function associated to $Z$ via  \eqref{zeta}.  By \cite[Theorem 1]{baetaludwigsemicontinuity} and \cite[Lemma~12]{baeta1d},  for any $u \in \LP(\mathbb{R}^n)$, the functional 
\begin{align*}
    Z_{\zeta}(u)= \int_{\mathbb{R}^n}\zeta(\det \D u(x))\dif x
\end{align*}
is a simple, $\tau$-upper semicontinuous, equi-affine and dually epi-translation invariant valuation that vanishes on indicator functions of polytopes and cylinder functions. Moreover, it satisfies
$$Z_{\zeta}(aq+\I_P)= \zeta(a)V_n(P)$$ 
for every $a\geq 0$ and  every polytope $P\subset \mathbb{R}^n$,  where $q(x)=\frac{1}{2}\sum_{i=1}^n x_i^2$.  Since both $Z$ and $Z_\zeta$ satisfy properties $(i)$--$(iv)$ of Proposition~\ref{uniquen}, it follows that $Z=Z_{\zeta}$. This completes the proof.
\end{proof}

To complete the proof of Theorem~\ref{equiv}, we proceed by induction on the dimension. The induction hypothesis is needed only in the proofs of the following two lemmas. More precisely, throughout these lemmas we assume that Theorem~\ref{equiv} holds in dimension $n-1$, in the following sense.

We say that Theorem~\ref{equiv} holds in dimension $n-1$ if, for every $\tau$-upper semicontinuous, equi-affine and dually epi-translation invariant valuation 
$\tilde{Z}:\LP(\mathbb{R}^{n-1})\rightarrow\mathbb{R}$, there exist  constants $d_0,d_1\in\mathbb{R}$ and a function $\gamma\in\conc$  such that
\begin{align}\label{(n-1)-dimensional}
   \tilde{Z}(u)=d_0+d_1\, V_{n-1}(\dom u)+\int_{\dom u} \gamma(\mbox{det}_{n-1}\D u(x_1,\dots, x_{n-1}))\dif x_1\dots \dif x_{n-1} 
\end{align}
 for every $u\in \LP(\mathbb{R}^{n-1})$.  Here, $(x_1,\dots, x_{n-1})$ are coordinates in $\mathbb{R}^{n-1}$, and $\dif x_1\dots \dif x_{n-1}$ denotes the $(n-1)$-dimensional Lebesgue measure. The term $\det_{n-1}\D u$ denotes the determinant of the Hessian matrix of $u$, where $u$ is defined on $\mathbb{R}^{n-1}$. For simplicity, we will often denote $(x_1,\dots, x_{n-1})$ and $\dif x_1\dots \dif x_{n-1}$ collectively by $x$ and $\dif x$, respectively.

\begin{lema}\label{simple}
Let $n\ge 2$ and assume that Theorem \ref{equiv} holds in dimension $n-1$. Let
$Z:\LP(\mathbb{R}^n)\rightarrow\mathbb{R}$
be a $\tau$-upper semicontinuous, equi-affine and dually epi-translation invariant valuation that vanishes on indicator functions of singletons. Then $Z$ is simple.
\end{lema}

\begin{proof}
Let $H\subset \mathbb{R}^n$ be a hyperplane. By equi-affine invariance of $Z$, we may assume without loss of generality that $H = \mathbb{R}^{n-1}$.
 Denote by $\mathcal{P}^{n-1}$ the set of polytopes contained in $\mathbb{R}^{n-1}$. Consider the restriction
$$\LP(\mathbb{R}^{n-1})=\{u\in\LP(\mathbb{R}^n)\mid \ \dom u \in \mathcal{P}^{n-1}\},$$
and define $\tilde{Z} = Z|_{\LP(\mathbb{R}^{n-1})}$. Then $\tilde{Z}$ is a $\tau$-upper semicontinuous, equi-affine and dually epi-translation invariant valuation on $\LP(\mathbb{R}^{n-1})$ which vanishes on indicator functions of singletons and is invariant under affine transformations of $\mathbb{R}^{n-1}$. Since Theorem~\ref{equiv} holds in dimension $n-1$ and $\tilde{Z}$ vanishes on indicator functions of singletons, we have
\begin{align}\label{eq10.9}
    Z(u) = \tilde{Z}(u) = d \,V_{n-1}(\dom u) + \int_{\dom u} \gamma(\mbox{det}_{n-1} \D u(x))  \dif x
\end{align}
for every $u \in \LP(\mathbb{R}^{n-1})$ with some $\gamma \in \conc$ and constant $d\in\mathbb{R}$.

Let $\varphi_t$ be the affine transformation which scales by $t>0$ in $\mathbb{R}^{n-1}$ and by $1/t^{n-1}$ in the orthogonal direction. By equi-affine invariance of $Z$,
\begin{align*}
    Z(u \circ \varphi_t) = Z(u).
\end{align*}
Applying \eqref{eq10.9} to $u \circ \varphi_t$, we obtain
\begin{align}
\nonumber Z(u) & =  d\, V_{n-1}(\dom u)+\int_{\dom u} \gamma(\mbox{det}_{n-1}\D u(x))\dif x\\
\nonumber& = d\,V_{n-1}(\varphi_t^{-1}(\dom u))+\int_{\varphi^{-1}_t(\dom u) } \gamma (\mbox{det}_{n-1} \D u\circ \varphi_t(x))\dif x \\
\nonumber & = \frac{d}{t^{n-1}}\, V_{n-1}(\dom u)+\int_{\varphi^{-1}_t(\dom u)} \gamma\left(t^{2(n-1)}\mbox{det}_{n-1} \D u( \varphi_t(x))\right)\dif x\\
\nonumber & = \frac{d}{t^{n-1}}\, V_{n-1}(\dom u)+\dfrac{1}{t^{n-1}}\int_{\dom u} \gamma \left(t^{2(n-1)}\mbox{det}_{n-1} \D u( x)\right)\dif x.
\end{align}
Hence
\begin{multline}\label{t}
\frac{d}{t^{n-1}}\, V_{n-1}(\dom u)+\dfrac{1}{t^{n-1}}\int_{\dom u} \gamma \left(t^{2(n-1)}\mbox{det}_{n-1} \D u( x)\right)\dif x\\= d\, V_{n-1}(\dom u)+ \int_{\dom u} \gamma(\mbox{det}_{n-1}\D u(x))\dif x    
\end{multline}
for every $t>0$ and $u\in \LP(\mathbb{R}^{n-1})$.

Taking $u = \I_P$, where $P \subset \mathbb{R}^{n-1}$ is a polytope with $V_{n-1}(P) > 0$, we obtain
$$\frac{d}{t^{n-1}}\, V_{n-1}(P)= d\, V_{n-1}(P)$$
for every $t>0$, since $\gamma(0)=0$. This immediately implies that $d = 0$.   

Now let $u(x) = \frac{1}{2}\|x\|^2+\I_P(x)$.  Then $\det_{n-1} \D u(x) =~1$ on $P$. With $d=0$,  \eqref{t} yields
$$
\frac{1}{t^{n-1}} \gamma(t^{2(n-1)})\, V_{n-1}(P)
= \gamma(1)\, V_{n-1}(P).
$$
Setting $\lambda = t^{2(n-1)}$, we obtain
$
\gamma(\lambda) = \gamma(1)\,\lambda^{1/2}.
$
Thus $\gamma \equiv c\sqrt{\cdot}$ for some  $c \ge 0$. 

Finally, consider the sequence 
$$u_k(x)=kx_1^2+\sum_{i=2}^{n-1}x_i^2+\I_{[0,1]^{n-1}\times\{0\}}(x), \quad k\in\mathbb{N},$$ 
which epi-converges to 
$$u(x)= \sum_{i=2}^{n-1}x_i^2+\I_{\{0\}\times[0,1]^{n-2}\times\{0\}}(x).$$

Since all functions are defined on $(n-1)$-dimensional domains, epi-convergence implies $\tau$-convergence in this setting without requiring uniform Lipschitz bounds. Hence $u_k$ is $\tau$-convergent to $u$, and we get
$$0\le \limsup_{k\rightarrow \infty} Z(u_k)= \limsup_{k\rightarrow \infty} c\, \sqrt{2^{n-1}k}\leq Z(u)=0,$$
and hence $c\equiv 0$. This completes the proof.
\end{proof}

\begin{lema}\label{cylinder}
Let $n\ge 2$ and assume that Theorem \ref{equiv} holds in dimension $n-1$. Let $Z: \LP(\mathbb{R}^n)\rightarrow\mathbb{R}$ be a simple, $\tau$-upper semicontinuous,  equi-affine and dually epi-translation invariant valuation that vanishes on indicator functions of polytopes in $\LP\RR$. Then $Z(u)= 0$ for every cylinder function $u\in\LP\RR$.
\end{lema}

\begin{proof}
Let $u \in \LP(\mathbb{R}^n)$ be a cylinder function.  
Thus, there exist $w \in \LP(\mathbb{R}^n)$, a closed interval $J\subset \mathbb{R}$,  
and an affine function $v$ such that
$$
u = w \,\square\, (v + \I_J),
$$
where $\dom w$ is a convex body of dimension at most $n-1$.

Since $Z$ is equi-affine and dually translation invariant, we may assume that $\dom w \subset H$, where $H$ is a hyperplane containing the origin,  
and that the interval $J$ lies in the line $L$.
We consider only segments $J$ lying in $L$. Set
$$\LP(H)=\bigl\{w\in\LP\RR\mid \ \dom w\subset \mathcal{P}(H)\bigr\}.$$  
Fix $w  \in \LP(H)$.  By the translation invariance of $Z$, the value of
$$
Z_1\big(w \,\square\, (v + \I_J)\big)
$$
remains unchanged if the closed line segment $J$ is translated along the   line,  since the epi-sum corresponds to the Minkowski sum of the epi-graphs. Consequently, $Z_1\big(w \,\square\, (v + \I_J)\big)$ depends only on the length of $J$, not on its position in $L$. Define
$$g_w(s)= Z_1(w\,\square\, (v+\I_{J})),\quad  s=V_1(J).$$

Let $J_1$ and $J_2$ be line segments in $L$, whose relative interiors are disjoint, and suppose that their  intersection consists of the origin  such that $v(0) = 0$. 
Set $s_1 = V_1(J_1)$ and $s_2 = V_1(J_2)$. Then
\begin{align*}
    (w\,\square\, (v+\I_{J_1}))\wedge (w\,\square\, (v+\I_{J_2})) & =w\,\square\, (v+\I_{J_1\cup J_2})\\
    (w\,\square\, (v+\I_{J_1}))\vee (w\,\square\, (v+\I_{J_2})) & =w.
\end{align*}
Since $Z_1$ is a translation invariant and simple valuation, we get
$$Z_1(w\,\square\, (v+\I_{J_1\cup J_2}))= Z_1(w\,\square\, (v+\I_{J_1}))+Z_1(w\,\square\, (v+\I_{J_2})),$$
which yields
$$g_w(s_1+s_2)=g_w(s_1)+g_w(s_2).$$
Thus, $g_w$ is an upper semicontinuous solution of Cauchy’s functional equation. Hence, there exists a real-valued function $w \mapsto Y(w)$ such that
$$g_w(s)= Y(w)\,s,$$
and therefore
$$Z_1(w\,\square\, (v+\I_{J}))=Y(w)V_1(J)$$
for every $w\in \LP(H)$ and every line segment $J$ in $L$. Note that $Y$ is defined on $\LP(H)$, and it is a $\tau$-upper semicontinuous,  equi-affine and  dually epi-translation invariant  valuation that vanishes on indicator functions of polytopes. Since Theorem \ref{equiv} holds in dimension $n-1$  and $Y$ vanishes on indicator functions of polytopes, we obtain
$$Y(w)= \int_P \gamma(\mbox{det}_{n-1} \D u(x))\dif x,$$
which implies
\begin{align}\label{v}
Z_1(w\,\square\, (v+\I_{J}))= V_1(J)\int_P \gamma(\mbox{det}_{n-1} \D u(x))\dif x, 
\end{align}
for every $u+\I_P\in \LP(H)$, every line segment $J$ in $L$, and some suitable $ \gamma\in\conc$.

As in the proof of Lemma \ref{simple},  choose equi-affine invariant transformations $\varphi_t$ that dilate
by a factor of $t>0$ in the direction of $H$ and by a factor of $1/t^{n-1}$ in the direction of $J$. Then
$$(w\,\square\, (v+\I_{J}))\circ \varphi_t$$
is a translate of the function $(u\circ \varphi_t+ \I_{\frac{1}{t}P})\,\square\, (v\circ\varphi_t+\I_{t^{n-1}J})$. From the equi-affine invariance of $Z_1$ and \eqref{v}, it follows that
\begin{align*}
    Z_1((u\circ \varphi_t+ \I_{\frac{1}{t}P})\,\square\, (v\circ\varphi_t+\I_{t^{n-1}J})) &= V_1(J)\int_{P} \gamma \left(t^{2(n-1)}\mbox{det}_{n-1} \D u( x)\right)\dif x\\
    &= Z_1(w\,\square\, (v+\I_{J}))\\
    &= V_1(J)\int_{P} \gamma \left(\mbox{det}_{n-1} \D u( x)\right)\dif x.   
\end{align*}
Hence $\gamma(t^{2(n-1)})= \gamma(t)$ for every $t>0$. Since $\gamma\in\conc$, we conclude that $\gamma\equiv 0$.
\end{proof}

Finally, combining Lemma~\ref{simple}, Lemma~\ref{cylinder}, and Proposition~\ref{prop}, we complete the proof of Theorem~\ref{equiv}. Recall that this theorem asserts that a functional
\(
Z:\LP\RR\rightarrow\mathbb{R}
\)
is a $\tau$-upper semicontinuous, equi-affine and dually epi-translation invariant valuation if and only if there exist constants $c_0,c_1\in\mathbb{R}$ and a function $\zeta\in\conc$ such that
\[
Z(u)=c_0+c_1\,V_n(\dom u)
+\int_{\dom u}\zeta(\det\D u(x))\,\dif x
\]
for every $u\in\LP\RR$.

\begin{proof}[Proof of Theorem \ref{equiv}]
By the discussion following Theorem~\ref{equiv}, it suffices to prove the necessity part of the theorem. We proceed by induction on the dimension $n$.

For  $n=1$, let $Z:\LP(\mathbb{R})\rightarrow\mathbb{R}$ be a $\tau$-upper semicontinuous, translation and  dually epi-translation   invariant valuation. By \cite[Theorem~2]{baeta1d}, there exist  constants $d_0,d_1$ and a  function $\gamma\in\conc$ such that
$$Z(u)=d_0+d_1\,V_1(\dom u)+\int_{\dom u}\gamma (u''(x))\dif x$$
for every $u\in\LP\NN$. Thus, the theorem holds for $n=1$.  

Assume now that the theorem holds in dimension $n-1$. We will show that it also holds in dimension $n$.

Let $x_0\in\mathbb{R}^n$. Since $Z$ is translation invariant, there exists a constant $c_0\in\mathbb{R}$ such that $Z(\I_{\{x_0\}})=c_0$ for every $x_0\in\mathbb{R}^n$. For $u\in \LP\RR$, define
\begin{align}\label{z_0}
    Z_0(u)=Z(u)-c_0.
\end{align}
Then $Z_0$ is a $\tau$-upper semicontinuous, equi-affine and  dually epi-translation  invariant valuation  on $\LP\RR$ that vanishes on indicator functions of singletons, i.e.,   $Z_0(\I_{\{x_0\}})=0$ for every $x_0\in\mathbb{R}^n$.  By Lemma~\ref{simple}, $Z_0$ is simple.

By Lemma \ref{non-positive},  there exists a constant $c_1\in\mathbb{R}$ such that
$Z_0(\I_P)= c_1\,V_n(P)$ for every polytope $P\subset\mathbb{R}^n$.
 Moreover, by \cite[Lemma 6 and Lemma 9]{baeta1d}, the functional 
$$
u \mapsto V_n(\dom u), \quad u\in \LP\RR,
$$ 
shares the same properties as $Z_0$, with the additional property of being $\tau$-continuous. Therefore,
\begin{align}\label{z1}
  Z_1(u)=Z_0(u)-c_1\,V_n(\dom u)  
\end{align}
defines a simple, $\tau$-upper semicontinuous,  equi-affine and dually epi-translation invariant valuation  that vanishes on indicator functions of polytopes in $\LP\RR$. By Lemma \ref{cylinder}, $Z_1$ also vanishes on every cylinder function in $\LP\RR$.

Since $Z_1:\LP\RR\to\mathbb{R}$ is a simple, $\tau$-upper semicontinuous, equi-affine and dually epi-translation invariant valuation that vanishes on indicator functions of polytopes and cylinder functions, Proposition~\ref{prop} yields a function $\zeta\in\conc$ such that
\begin{align}\label{z3}
    Z_1(u)= \int_{\dom u}\zeta(\det \D u(x))\dif x
\end{align}
for every $u\in \LP\RR$.

Finally, combining \eqref{z_0}, \eqref{z1}, and \eqref{z3}, we obtain
\begin{align*}
Z(u)
&=c_0+Z_0(u)\\
&=c_0+c_1\,V_n(\dom u)+Z_1(u)\\
&=c_0+c_1\,V_n(\dom u)+\int_{\dom u}\zeta(\det\D u(x))\,\dif x
\end{align*}
for every $u\in\LP\RR$. This completes the proof.
\end{proof}

\section*{Acknowledgments}
The author is grateful to Monika Ludwig for her invaluable help in discussions, for her careful reading of the manuscript, and for her insightful suggestions. The author thanks the referee for their careful reading of the manuscript and for helpful comments. This project was supported, in part, by the Austrian Science Fund (FWF) Grant-DOI: 10.55776/P34446 and Grant-DOI: 10.55776/P37030. For open access purposes, the author has applied a CC BY public copyright license to any author-accepted manuscript version arising from this submission.

\bibliography{gc}
\bibliographystyle{amsplain}

\newpage

\appendix

\section{Appendix}\label{apendice}

The proof of Proposition~\ref{prop10.3} relies on several auxiliary results whose proofs are technically involved and are therefore deferred to this appendix in order to improve the readability of the main text.

We establish Proposition~\ref{regular}, then Proposition~\ref{regular2}, and finally Proposition~\ref{singular}.

\subsection{Proof of Proposition \ref{regular}}

We recall the statement of Proposition~\ref{regular}. Moreover, $N \subset \dom u$ denotes the set of interior points at which $u$ is twice differentiable.

\begin{prop}\label{regularr}
Let $Z:\LP\RR\rightarrow\mathbb{R}$ be a simple, $\tau$-upper semicontinuous, equi-affine and  dually epi-translation   invariant valuation that   vanishes on  indicator functions of polytopes and cylinder functions. Let $u \in \LP(\mathbb{R}^n)$,
let $x_0\in N$ be such that $\det\D u(x_0)>0$,  let $\rho>0$, and let $l^i, l^c$  be the functions satisfying \eqref{li_lc}. Then there exist $r(x_0)>0$, 
polytopes $P_r(x_0)$ containing $x_0$ in their interiors, and  functions $v_r^{x_0}\in P\LQ\RR$ such that
\begin{align}\label{eq10.47a}
  Z(u+\I_{P_r(x_0)}) \leq Z(v_r^{x_0}+\I_{P_r(x_0)}) + \frac{\rho}{2}V_n(P_r(x_0))
\end{align}
and
\begin{align}\label{eq10.48b}
l^c(x)\leq v_r^{x_0}(x)\leq l^i(x)    
\end{align}
for all $x\in \mathbb{R}^n$ and $0<r\leq r(x_0)$. Moreover, for $0<r\leq r(x_0)$,  the function $v_r^{x_0}$ can be chosen so such that
$u$ and $v_r^{x_0}$ share a common supporting hyperplane at $x_0$, and 
\begin{align}\label{eq4111c}
 \{x\in \mathbb{R}^n\mid v_r^{x_0}(x)\geq u(x)\}   \subset P_r(x_0), 
\end{align}
while the  polytopes $P_r(x_0)$ satisfy 
\begin{align}\label{eq10.26..d}
\diam(P_r(x_0))^n= n^{\frac{n}{2}}\, V_n(P_r(x_0)).    
\end{align}
\end{prop}

Throughout this section, we assume that
$Z:\LP\RR\rightarrow\mathbb{R}$ is a simple, $\tau$-upper semicontinuous, equi-affine and dually epi-translation invariant valuation that vanishes on indicator functions of polytopes and cylinder functions.

Let $x_0\in N$ be a point such that $\det \D u(x_0)>0$ and let $\rho>0$. Because of the equi-affine and vertical translation invariance of $Z$, we may assume without loss of generality that
$$
x_0 = 0, \qquad u(0) = 0, \qquad \nabla u(0) = 0, \qquad \D u(0) = \lambda \diag(1, \dots, 1),
$$
so that all eigenvalues of $\D u(0)$ are equal to $\lambda > 0$. Let $t \in (0,1)$. 
Consider quadratic functions $q_{u,t}^i$ and $q_{u,t}^c$ defined by
$$
q_{u,t}^i(x) = (1+t)\lambda\, q(x), 
\qquad
 q_{u,t}^c = (1-t) \lambda\, q(x),
$$
where $q$ is given by \eqref{q}. Then
$$
q_{u,t}^i(x_0) = u(x_0) = q_{x_0,t}^c(x_0),
$$
and both $q_{u,t}^i$ and $q_{u,t}^c$ share the same tangent hyperplane as $u$ at $x_0$.  In a sufficiently small neighborhood $V(x_0)$ of $x_0$, we have
$$
(q_{u,t}^c  + \I_{V(x_0)})(x) 
\le (u + \I_{V(x_0)})(x) 
\le (q_{u,t}^i + \I_{V(x_0)})(x).
$$
For simplicity of notation, we  denote $q_{u,t}^i$ by $q_u$.

Let $t>0$ be sufficiently small such that, if $\zeta \not\equiv 0$, we have 
\begin{align}\label{10.24}
  4\sqrt{t}< 1 - \max\left\{0, \left(1-\frac{\rho}{4\zeta(2^n\lambda^n)}\right)\right\}^{1/n} 
\end{align}
and 
\begin{align}\label{10.25}
    2\sqrt{t}<\left(1+\frac{\rho}{4\zeta(2^n\lambda^n)}\right)^{1/n}-1.
\end{align}
Note that $\rho>0$ is fixed and does not depend on $x_0$, and that, by \eqref{10.24}, $0<t<\frac{1}{16}$.

Recall that $C=[-1,1]^n$. Moreover, 
\begin{align}\label{eq10.26.}
   \diam(C)^n= n^{\frac{n}{2}}\, V_n(C).
\end{align}

Note that for $r>0$ sufficiently small the polytope $rC$ is contained in $\dom u$ and, by construction, 
$$
(q_{u,t}^c + \I_{rC})(x) \leq (u + \I_{rC})(x) \leq (q_u + \I_{rC})(x) \quad \text{for all } x \in \mathbb{R}^n.
$$

We aim to show that $Z(q_u+\I_{rC})$ is nearly as large as $Z(u + \I_{rC})$ for  sufficiently small $r>0$.  Consider the convex hull
$$\conv(\epi (u + \I_{rC})\cup\epi  (q_u+\I_C)).$$ 
By \eqref{set_func}, there exists a function  in $\LP(\mathbb{R}^n)$ corresponding to this convex hull,  which we will denote by
$$
\conv(u + \I_{rC}, q_u + \I_C),
$$
and which differs from $q_u + \I_C$ only near $x_0$.  We will construct a convex function  $w_r$ using the part of $\conv(u + \I_{rC}, q_u + \I_C)$, where it does not coincide with $q_u + \I_C$,   in such a way that $w_r \seq q_u+\I_C$ as $r\rightarrow 0$. Then, by the $\tau$-upper semicontinuity of $Z$,  for sufficiently small $r>0$,   $Z(w_r)$ is close to $Z(q_u+\I_C)$, and the construction ensures that $Z(u + \I_{rC})$ is also close to $Z(q_u + \I_{rC})$.  Finally, we will define $v_r^{x_0}\in P\LQ\RR$ by taking $(q_u+\I_{rC})(x)$ for $x\in rC$ and using affine and cylinder functions for $x\notin rC$.   Having described the plan, we now proceed with the construction of $w_r$.

Throughout the text, we will use the notation $[p_1,p_2]$ to denote the closed line segment with endpoints $p_1$ and $p_2$.

\begin{lema}\label{claim1}
For $0 < r < 1$, define
$$
G_r = \{ y \in \mathbb{R}^n \mid \conv(u + \I_{rC}, q_u + \I_C)(y) \neq (q_u + \I_C)(y) \}.
$$
Then
\begin{align}\label{eq10.29}
G_r \subset (1 + 2\sqrt{t})\, rC.
\end{align}    
\end{lema}

\begin{proof}
Let $p_1=(y,y_{n+1})$, $y_{n+1}\neq q_u(y)$, and  $p_2= (\Bar{y},q_u(\Bar{y}))$  such that the segment $[p_1,p_2]$ lies in the tangent hyperplane to $q_u$ at $\bar{y}$. Write
$\Bar{y}=(\Bar{y}_1, \dots, \Bar{y}_{n})$ and  $y=(y_1,\dots, y_n)$. Then
\begin{align*}
    y_{n+1}-\lambda(1+t)(\Bar{y}_1^2+\cdots + \Bar{y}_n^2)= 2\lambda(1+t) \langle (\Bar{y}_1, \dots, \Bar{y}_n), (y_1-\Bar{y}_1, \dots, y_n-\Bar{y}_n)\rangle,
\end{align*}
which implies
\begin{align*}
   \|y-\Bar{y}\|^2
    & = \|y\|^2-\frac{y_{n+1}}{\lambda(1+t)},
\end{align*}
and for $y_{n+1}=q_{u,t}^c(y)$,
\begin{align}\label{eq10.32}
 \|y-\bar{y}\| = \sqrt{\frac{2 t\|y\|^2}{1+t}}\le \sqrt{2t}\|y\|.  
\end{align}

Hence, we have the inclusion
\begin{align}\label{10.31}
\bigl\{ \conv(q_{u,t}^c(y)+\I_{\{(y,q_{u,t}^c(y))\}}, q_u+\I_C)\neq q_u+\I_C\bigr\}\subset y+\dfrac{\sqrt{2 t}}{\sqrt{n}}\|y\|\, C.
\end{align}
Since $q_{u,t}^c(y)\leq (u + \I_{rC})(y)$, combining \eqref{eq10.32} and \eqref{10.31} yields
$$\bigl\{\conv((u + \I_{rC})(y)+\I_{\{(y,(u + \I_{rC})(y)\}}, q_u+\I_C)\neq q_u+\I_C\bigl\}\subset 
y+\sqrt{t}(1+2\sqrt{t})\,rC$$
for every $y\in rC$. Moreover, since $y \in  rC$ and by \eqref{10.24} we have $t < \frac{1}{16}$, this completes the proof of the lemma.
\end{proof}

\begin{lema}\label{lema10.1}
Under the above hypotheses,
$$Z(u + \I_{rC})\leq Z(q_u+\I_{rC})+\frac{\rho}{4}V_n(rC)$$
for all sufficiently small $r > 0$.
\end{lema}

\begin{proof}
Assume, by contradiction, that for arbitrarily small $r>0$,
\begin{align}\label{eq10.28}
Z(u + \I_{rC})> Z(q_u+\I_{rC})+\frac{\rho}{4}V_n(rC).
\end{align}

By Lemma~\ref{claim1}, for sufficiently small $r>0$, the set
\[
G_r=\{y\in\mathbb{R}^n:\conv(u+\I_{rC},q_u+\I_C)(y)\neq(q_u+\I_C)(y)\}
\]
is contained in $(1+2\sqrt{t})rC$. Consequently,
$\conv(u+\I_{rC},q_u+\I_C)$ differs from $q_u+\I_C$ only in the neighborhood
$(1+2\sqrt{t})rC$ of $x_0$.

Let $m_r$ be the maximal number of points $y_1,\dots,y_{m_r}\in C$ for which the translates
\begin{align}\label{eq10.33}
  y_i+(1+2\sqrt{t})r\,C
\end{align}
are pairwise disjoint and contained in $C$. In particular,
$$
m_r \le \frac{V_n(C)}{V_n\big((1+2\sqrt{t})rC\big)}.
$$
Then
\begin{align}\label{eq10.34}
   m_rV_n((1+2\sqrt{t})r\, C)\rightarrow V_n(C) \quad \text{as } r\rightarrow 0.
\end{align}
 
Let $z\in\mathbb{R}^n$  and consider a quadratic function 
$$
\bar{q}(x) = \sum_{i=1}^n \gamma_i x_i^2 + \sum_{i=1}^n \beta_i x_i + \alpha,
$$  
where $\gamma_i\geq 0$ for   $i = 1, \dots, n$, $\beta_i \in \mathbb{R}$, and $\alpha \in \mathbb{R}$.  Define the function $\phi:\mathbb{R}^n\times\mathbb{R}\rightarrow \mathbb{R}^n\times\mathbb{R}$  by
\begin{align}\label{affine}
  \phi(x,y) = (z+x, y+\sum_{i=1}^n\gamma_iz_i^2+\sum_{i=1}^n\beta_i z_i+2\sum_{i=1}^n\gamma_iz_ix_i) =(z+x,y+\ell_z(x)),
\end{align}
where 
\begin{align}\label{psi}
 \ell_z(x)= 2\sum_{i=1}^n\gamma_iz_ix_i +\sum_{i=1}^n\gamma_iz_i^2+\sum_{i=1}^n\beta_i z_i   
\end{align} 
is  affine in $x$. Note that $\phi(x,\bar{q}(x))=(z+x,\bar{q}(z+x))$ defines a $C^1$ map.

For each $i = 1, \dots, m_r$, define
$$g_i(x+y_i)=(u + \I_{rC}+ \ell_{y_i})(x),$$
where $\ell_{y_i}$ is the affine function associated with $y_i$ and $q_u$, as in \eqref{psi}.  Next, define
$$w_r=\conv(g_1, \dots, g_{m_r},q_u+\I_C),$$
where the convex hull is taken in terms of the epi-graphs, i.e.,
$$
\epi(w_r) = \conv\Bigl(\epi(g_1) \cup \dots \cup \epi(g_{m_r}) \cup \epi(q_u + \I_C)\Bigr).
$$ 
By construction, 
$$w_r\seq q_u+\I_{C}$$
as $r\rightarrow 0$, and  due to $\eqref{eq10.29}$ and \eqref{eq10.33},  
\begin{align}\label{eq10.35}
w_r(x+y_i)\leq (u + \I_{rC}+\ell_{y_i})(x)
\end{align}
 for $i=1,\dots, m_r, x\in\mathbb{R}^n$, and  sufficiently small $r>0$.

Next, we dissect
$$
C \setminus \bigcup_{i=1}^{m_r} \bigl( y_i +  rC \bigr)
$$
into polytopes $P_1, \dots, P_{k_r}$ whose interiors are disjoint from each $y_i +  rC$, $i=1,\dots,m_r$. It follows from \eqref{eq10.35} that the intersections $(y_i + rC) \cap (y_j + rC)$ are empty for $i \neq j$.

Since $Z$ is a simple, non-negative, $\tau$-upper semicontinuous, and dually epi-translation invariant valuation,  we have $\zeta((1+t)^n\lambda^n)\geq 0$.    Using formulas \eqref{mudanca_dominio}, \eqref{mudanca_zeta}, and  by \eqref{eq10.28}, \eqref{eq10.29}, we obtain, for every $\eta > 0$,
\begin{align*}
    Z(w_r) &=  \sum_{i=1}^{m_r}Z(u + \I_{rC}+\ell_{y_i}) +  \sum_{j=1}^{k_r}Z(w_r+\I_{P_j})\\
    &\geq  m_r Z(u + \I_{rC}) + \dfrac{Z(q_u+\I_C)}{V_n(C)}(V_n(C)-m_rV_n((1+2\sqrt{t})\,rC))\\
    &\geq Z(q_u+\I_C)+m_r\left(Z(q_u+\I_{rC})+\frac{\rho}{4}V_n(rC)\right) -\zeta((1+t)^n\lambda^n)m_r(1+2\sqrt{t})^nV_n(rC)\\
    &\geq Z(q_u+\I_C) + m_rV_n(rC)\left(\dfrac{\rho}{4}-\zeta((1+t)^n\lambda^n)(1+2\sqrt{t})^n+\zeta((1+t)^n\lambda^n)\right).
\end{align*}

By \eqref{eq10.29}, we have $m_rV_n(rC)\leq m_rV_n(G_r)$, and from \eqref{eq10.34} it follows that
$$m_rV_n(rC)= m_rV_n(rC)\rightarrow  \dfrac{V_n(C)}{(1+2\sqrt{t})^n}$$
as $r\rightarrow 0$. Moreover, by \eqref{10.25}, 
$$2\sqrt{t}< \left(1+\dfrac{\rho}{4\zeta(2^n\lambda^n)}\right)^{1/n}-1,$$ 
which immediately gives  
$$ \dfrac{\rho}{4}-\zeta((1+t)^n\lambda^n)(1+2\sqrt{t})^n+\zeta((1+t)^n\lambda^n)> 0,$$
because $t<1$ and $\zeta$ is non-decreasing. 
Since $Z$ is $\tau$-upper semicontinuous and $w_r$ is  $\tau$-convergent to $q_u+\I_C$ as $r\rightarrow 0$, we obtain
\begin{align*}
 Z(q_u+\I_C)   & \geq \limsup_{r\rightarrow 0} Z(w_r)\\
 & \geq Z(q_u+\I_C)  + \limsup_{r\rightarrow 0}   m_rV_n(rC)\left(\dfrac{\rho}{4}-\zeta((1+t)^n\lambda^n)(1+2\sqrt{t})^n+\zeta((1+t)^n\lambda^n)\right)\\
 & \geq  Z(q_u+\I_C) +\dfrac{V_n(C)}{(1+2\sqrt{t})^n}\left(\dfrac{\rho}{4}-\zeta((1+t)^n\lambda^n)(1+2\sqrt{t})^n+\zeta((1+t)^n\lambda^n)\right),
\end{align*}
which is a contradiction. Therefore, we conclude that
\begin{align}\label{eq10.38}
Z(u + \I_{rC})\leq Z(q_u+\I_{rC})+\frac{\rho}{4}V_n(rC)
\end{align}
holds for  sufficiently small $r>0$.
\end{proof}

The goal is now to construct, for small $r>0$, a function $v_r\in P\LQ(\mathbb{R}^n)$ that coincides with $q_u+\I_{rC}$ on a smaller subset of $rC$, is composed outside this subset of affine and cylinder pieces, and agrees with $u$ at a prescribed point of $rC$.

Let $K \subset \mathbb{R}^n$ be a convex body and let $\tilde{q}:\mathbb{R}^n \to (-\infty, \infty]$ be a piecewise linear-quadratic  function such that $K\subset \interior (\dom \tilde{q})$. We define $\mathrm{ext}_K(\tilde{q}): \mathbb{R}^n \to \mathbb{R}$ by 
\begin{align}\label{tang}
\mathrm{ext}_K(\tilde{q})(y)=
\begin{cases}
\tilde{q}(y), & y \in K,\\
\displaystyle \sup_{x \in \bd K} \left\{ \tilde{q}(x) + \langle \nabla \tilde{q}(x), y - x\rangle \right\}, & y \notin K.
\end{cases}
\end{align}
Since $\tilde{q}$ is of class $C^1$ in $\interior (\dom \tilde{q})$, for each $x \in \bd K$ the affine function
\begin{align*}
\ell_x(y) = \tilde{q}(x) + \langle \nabla \tilde{q}(x), y - x\rangle
\end{align*}
defines the supporting hyperplane to $\tilde{q}$ at $x$. Moreover, for every $y \notin K$, the supremum is attained. Indeed, the map
\begin{align*}
x \mapsto \tilde{q}(x) + \langle \nabla \tilde{q}(x), y - x\rangle
\end{align*}
is continuous on $\bd  K$, which is compact. Therefore, there exists $x \in \bd K$ such that
\begin{align*}
\mathrm{ext}_K(\tilde{q})(y)= \tilde{q}(x) + \langle \nabla \tilde{q}(x), y - x\rangle.
\end{align*}

\begin{lema}\label{claim2}
For every $r$, $0 < r < 1$, define
$$
G_r = \{ y \in \mathbb{R}^n \mid \mathrm{ext}_{(1-4\sqrt{t}) rC}(q_u)(y) = u (y) \}.
$$
Then
\begin{align*}
G_r \subset rC.
\end{align*}    
\end{lema}

\begin{proof}
We proceed as in   Lemma \ref{claim1}. Let  $\bar{y}\in \bd((1-4\sqrt{t})\,rC)$,  and let $y\in \dom u$ be such that the segment $[(y, u(y)),(\Bar{y}, q_u(\Bar{y}))]$ lies in the tangent hyperplane to $q_u$ at $\Bar{y}$. Then by   \eqref{eq10.32},
$$\|y-\Bar{y}\|= \sqrt{t}\|y\|.$$
By assumption, $\|\bar{y}\|\leq (1-4\sqrt{t})r\sqrt{n}$, and the triangle inequality gives
$$\|y\|\leq \|y- \Bar{y}\|+\|\Bar{y}\|\leq \sqrt{t}\|y\|+(1-4\sqrt{t})r\sqrt{n}.$$
Rearranging, we obtain
$$\|y\|\leq \dfrac{(1-4\sqrt{t})r\sqrt{n}}{1-\sqrt{t}},$$
and consequently
$$y\in \bar{y}+\sqrt{t}\dfrac{1-4\sqrt{t}}{1-\sqrt{t}}\sqrt{n}\dfrac{1+2\sqrt{t}}{\sqrt{n}}rC.$$
Recalling that $\bar{y}\in (1-4\sqrt{t})\,rC$, and using the fact that $8t-2\sqrt{t}+3>0$ for all $t\geq 0$, this completes the proof of the claim.
\end{proof}

\begin{lema}\label{lema*}
For sufficiently small $r>0$, there exists  $v_r\in P_{\LQ}\RR$ such that
\begin{align*}
 (u + \I_{rC})(x)\leq (v_r+\I_{rC})(x) \leq  (q_u+\I_{rC})(x) 
\end{align*}
and 
\begin{align}\label{eq10.39}
    Z(u + \I_{rC})\leq Z(v_r+\I_{rC})+\dfrac{\rho}{2}V_n(rC).
\end{align}
\end{lema}

\begin{proof}
By Lemma~\ref{claim2}, for sufficiently small $r>0$, the tangential extension of $q_u+ \I_{(1-4\sqrt{t})\,rC}$ (see \eqref{tang}) coincides with $u$ on a subset of $rC$.

Now define 
\begin{align}\label{eq10.42}
v_r = \mathrm{ext}_{(1-4\sqrt{t}) rC}(q_u) \vee  l^c, 
\end{align}
where $l^c$ is the piecewise affine function chosen in \eqref{li_lc}. Note that $v_r \in P\LQ\RR$ and since $Z$ vanishes on cylinder functions, is dually epi-translation invariant and vanishes on indicator functions of polytopes,
\begin{align}\label{eq10.43}
 Z(v_r) = Z(q_u + \I_{(1-4\sqrt{t})\,rC}).  
\end{align}
By \eqref{eq10.38}, \eqref{mudanca_dominio}, \eqref{zeta} and using that $Z$ is a simple valuation, we get
\begin{align}
    \nonumber  Z(u + \I_{rC}) & \leq Z(q_u+\I_{rC})+\frac{\rho}{4}V_n(rC)\\
    \nonumber & = Z(q_u+\I_{(1-4\sqrt{t})\,rC})\frac{V_n(rC)}{V_n((1-4\sqrt{t})\,rC)}+\frac{\rho}{4}V_n(rC)\\
    \nonumber & = Z(q_u+\I_{(1-4\sqrt{t})\,rC}) +\frac{\rho}{2}V_n(rC)\\
     & \qquad + V_n(rC)\left(\zeta((1+t)^n\lambda^n) - \frac{\rho}{4} - \zeta((1+t)^n\lambda^n)(1-4\sqrt{t})^n \right).\label{eq**}
\end{align}

By \eqref{10.24}, if $\rho<4\zeta(2^n\lambda^n)$,  we have
$$
4\sqrt{t}< 1 - \left(1-\frac{\rho}{4\zeta(2^n\lambda^n)}\right)^{1/n},
$$
which together with $t<\frac{1}{16}$ implies that
\begin{align}\label{neg}
\zeta((1+t)^n\lambda^n)-\frac{\rho}{4}-\zeta((1+t)^n\lambda^n)(1-4\sqrt{t})^n < 0.  
\end{align}
On the other hand, if $\rho \ge 4\zeta(2^n\lambda^n)$, then $1 - 4\sqrt{t} > 0$, which again implies 
\begin{align}\label{neg+}
\zeta((1+t)^n\lambda^n)-\frac{\rho}{4}-\zeta((1+t)^n\lambda^n)(1-4\sqrt{t})^n < 0.  
\end{align}

Combining \eqref{eq10.43}, \eqref{eq**}, \eqref{neg}, and \eqref{neg+}, we obtain
$$
Z(u + \I_{rC}) \leq Z(v_r) + \frac{\rho}{2}V_n(rC),
$$
as desired.
\end{proof}

We are now ready to prove Proposition~\ref{regularr}.

\begin{proof}[Proof of Proposition~\ref{regularr}]
Let $r>0$ be sufficiently small so that the conclusions of Lemma~\ref{lema10.1} and Lemma~\ref{lema*} hold. Define
\begin{align}\label{eq10.44}
    P_r(x_0)=rC
    \quad \text{and} \quad
    v_r^{x_0}=v_r.
\end{align}
By Lemma~\ref{lema*}, we have $v_r^{x_0}\in P\LQ\RR$ and
\[
Z(u+\I_{P_r(x_0)})
\le
Z(v_r^{x_0}+\I_{P_r(x_0)})
+\frac{\rho}{2}V_n(P_r(x_0)).
\]
This gives \eqref{eq10.47a}.

Moreover, by the construction of $v_r$ in \eqref{eq10.42}, and by choosing $r>0$ sufficiently small, we have
\[
l^c(x)\leq v_r^{x_0}(x)\leq l^i(x)
\quad \text{for all } x\in\mathbb{R}^n,
\]
which gives \eqref{eq10.48b}. Also, since $v_r$ is defined through the tangential extension of $q_u$ on $(1-4\sqrt{t})rC$ and shares the same supporting hyperplane as $u$ at $x_0$, the functions $u$ and $v_r^{x_0}$ share a common supporting hyperplane at $x_0$. By Lemma~\ref{claim2}, the set where this extension can coincide with $u$ is contained in $rC=P_r(x_0)$; hence
\[
\{x\in\mathbb{R}^n: v_r^{x_0}(x)\geq u(x)\}
\subset P_r(x_0),
\]
which gives \eqref{eq4111c}.

Finally, by \eqref{eq10.26.}, since $P_r(x_0)=rC$, we have
\[
\diam(P_r(x_0))^n
=
n^{\frac n2}V_n(P_r(x_0)).
\]
Thus there exists $r(x_0)>0$ such that \eqref{eq10.47a}, \eqref{eq10.48b}, \eqref{eq4111c}, and \eqref{eq10.26..d} hold for every $0<r\le r(x_0)$. This completes the proof.
\end{proof}

\subsection{Proof of Proposition \ref{regular2}}
We next recall the statement of Proposition~\ref{regular2}.

\begin{prop}\label{regularr2}
Let $Z:\LP\RR\rightarrow\mathbb{R}$ be a simple, $\tau$-upper semicontinuous, equi-affine and  dually epi-translation   invariant valuation that   vanishes on  indicator functions of polytopes and cylinder functions. Let $u \in \LP(\mathbb{R}^n)$,
let $x_0\in N$ be such that $\det \D u(x_0)=0$,  let $\rho>0$, and let $l^i, l^c$  be the functions satisfying \eqref{li_lc}. Then there exist $r(x_0)>0$ and   polytopes $P_r(x_0)$ containing $x_0$ in their interiors  such that
\begin{align}\label{eq10.67m}
  Z(u+\I_{P_r(x_0)})\leq \frac{\rho}{2}\, V_n(P_r(x_0))
\end{align}
and
\begin{align}\label{eq10.66m}
 l^c(x)\leq v_r^{x_0}(x)\leq l^i(x)    
\end{align}
for all  $x\in \mathbb{R}^n$ and $0<r\leq r(x_0)$, where
\begin{align}\label{eq41112m}
  v_r^{x_0}= l_{x_0}\vee l^c,  
\end{align}
and $l_{x_0}$ is the affine function defining the tangent hyperplane to $u$ at $x_0$, coinciding with $u$ at $x_0$. Moreover, for $0<r\leq r(x_0)$, the 
polytopes $P_r(x_0)$ satisfy 
\begin{align}\label{eq10.53..m}
\diam(P_r(x_0))^n= n^{\frac{n}{2}}\, V_n(P_r(x_0)).    
\end{align}
\end{prop}

Let $x_0\in N$ be a point such that $\det \D u(x_0)=0$ and $\rho>0$. As in the proof of Proposition \ref{regularr}, since $Z$ is  equi-affine and vertical translation invariant, we may assume without loss of generality that
\begin{align}\label{assump}
 x_0 = 0, \qquad u(0) = 0, \qquad \nabla u(0) = 0.   
\end{align}
Since $u$  is twice differentiable at $x_0$, by Taylor's expansion we can write, for $x$ near $x_0$,
\begin{align*}
    u(x)=\frac{1}{2}\langle Ax, x\rangle + o(\|x\|^2),
\end{align*}
where $A$ is a  positive semidefinite   matrix and $\frac{o(\|x\|^2)}{\|x\|^2}\rightarrow 0$ as $\|x\|\rightarrow 0$.
Since $A$ is symmetric and positive semidefinite, there exist an orthogonal matrix $U$ and a diagonal matrix  $D=~\diag(\lambda_1, \dots, \lambda_n)$ with $\lambda_i\geq 0$ for all $i=1,\dots, n$, such that $A= U DU^{-1}$, and this implies that
$$\langle A x, x\rangle = \langle Dx, x\rangle= \sum_{i=1}^n \lambda_ix^2.$$
Since $\det A = 0$, we have $\det D = 0$, which means that there exists $k\geq 1$ such that, without loss of generality, $\lambda_1= \cdots = \lambda_k=0$ and $\lambda_{k+1}, \dots, \lambda_n>0$. Consequently, there is an equi-affine map $\varphi$ such that  
$$
\lambda_1( u\circ \varphi^{-1}(\varphi(x_0)))= \cdots =\lambda_k( u\circ \varphi^{-1}(\varphi(x_0)))=0
$$
and 
$$
\lambda_{k+1}( u\circ \varphi^{-1}(\varphi(x_0)))= \dots = \lambda_{n}( u\circ \varphi^{-1}(\varphi(x_0)))=1.
$$
Therefore, without loss of generality, due to the $\V\N$ invariance of $Z$, we may assume that  $A=\diag(0,\dots, 0, 1, \dots, 1)$, where the entry $0$ appears $k$ times.

Recall that 
\begin{align}\label{eq10.53.}
   \diam(C)^n= n^{\frac{n}{2}}\, V_n(C).
\end{align}

Next, consider the quadratic function 
\begin{align*}
   q_t(x)= \dfrac{1}{2}\sum_{i=1}^k tx_i^2+\dfrac{1}{2}\sum_{j=k+1}^n(1+t)x_j^2 
\end{align*}
for $t>0$. Since $Z$ is  non-negative, simple,  and  $\tau$-upper semicontinuous, and since $q_t$ is $\tau$-convergent  to $\frac{1}{2}\sum_{j=k+1}^nx_j^2$ as $t\rightarrow 0$, we obtain
\begin{align*}
    \limsup_{t\rightarrow 0} Z(q_t+\I_C)=0.
\end{align*}
Hence, for sufficiently small $t$,
\begin{align}\label{eq10.52}
Z(q_t+\I_C)\le \dfrac{\rho \, V_n(C)}{4(4n)^n}.    
\end{align}
We assume that $t$ is chosen to satisfy \eqref{eq10.52}  and $0 < t < 1$. Our goal is to show that $Z(u+\I_{rC})$ is not much larger than $Z(q_t+\I_{rC})$ for $r>0$ sufficiently small.

\begin{lema}\label{claim3}
Let
\[
\bar{q}_t(x)=\frac{t}{2}\sum_{i=1}^k x_i^2
+\frac{1+t}{2}\sum_{j=k+1}^n x_j^2.
\]
There exists $r_1>0$ such that for $0<r\le r_1$,
$
\{u+\I_{rC} = \bar{q}_t\} = \{x_0\}
$.
Define
$$
G_r = \{y\in\mathbb{R}^n\mid  \conv(u+\I_{rC}, \bar{q}_t+\I_C)(y) \neq (\bar{q}_t+\I_C) (y)\}.
$$
Then
$$
G_r \subset  4 n rC.
$$    
\end{lema}

\begin{proof}
First, note that
\begin{align}\label{eq10.56}
  \dfrac{t}{4} \sum_{i=1}^k x_i^2 + \dfrac{1}{2} \sum_{j=k+1}^nx_j^2 \le q_t(x) \le \dfrac{3t}{2} \sum_{i=1}^k x_i^2 + \dfrac{(1+t\, n)}{2} \sum_{j=k+1}^nx_j^2.
\end{align}

Let $(y,0)$ be a point in $\mathbb{R}^n=H(x_0,(\bar{q}_t+ \I_C)(x_0))$, the 
 hyperplane tangent to $\bar{q}_t$ at $x_0$. Denote by
$$\dist((y,0), \epi (\bar{q}_t+ \I_C))= \inf\bigl\{\|(y,0)-(x,(\bar{q}_t+ \I_C)(x))\|:x\in C\bigr\}$$ 
the distance  from $(y,0)$ to the epi-graph of $\bar{q}_t+\I_C$.
Since $\bar{q}_t$ is  convex,
$$\dist((y,0), \epi (\bar{q}_t+ \I_C))\geq \dist((y,0), H(y,\bar{q}_t(y))).$$
As the angle between $\mathbb{R}^n$ and $H(y,\bar{q}_t(y))$ tends to $0$ as $\|y\|\to 0$, it follows that
$$\dist((y,0), H(y,\bar{q}_t(y))) \geq \frac{1}{2}\bar{q}_t(y)$$
for $\|y\|$ sufficiently small. Hence,
\begin{align}\label{eq10.57}
 \dist((y,0), H(y,\bar{q}_t(y)))  \geq    \dfrac{t}{8} \sum_{i=1}^k y_i^2 + \dfrac{1}{4} \sum_{j=k+1}^ny_j^2
 \geq \dfrac{t}{8}\|y\|^2
\end{align}
holds for $\|y\|$ sufficiently small.

By \eqref{assump},  we can represent $u$ in a neighborhood of the origin as
$$u(x) = \dfrac{1}{2}\sum_{j=k+1}^nx_j^2+ o(\|x\|^2),$$
where $\frac{o(\|x\|^2)}{\|x\|^2}\to 0$ as $\|x\|\to 0$. By convexity of $u$, there exists $r_1>0$ such that
\begin{align}\label{eq10.58}
    u(x)\geq \frac{(1-t)}{2}\sum_{j=k+1}^nx_j^2 
\end{align}
for $\|x\|\leq r_1$.

As in the proof of Claim 1 in Lemma \ref{lema10.1}, let $y\in rC$ and let $(\Bar{y}, \bar{q}_t(\Bar{y}))$ be such that $[(y,(u+\I_{rC})(y)), (\Bar{y}, \bar{q}_t(\Bar{y}))]$ lies in the hyperplane tangent to $\bar{q}_t$ at $\Bar{y}$. By \eqref{eq10.56}, \eqref{eq10.58}, and the inequality $\dist((y,(u+\I_{rC})(y)), \epi (\bar{q}_t+ \I_C))\leq \bar{q}_t(y)-(u+\I_{rC})(y)$, we obtain 
\begin{align}
 \nonumber\dist((y,(u+\I_{rC})(y)), \epi (\bar{q}_t+ \I_C)) & \leq    \dfrac{3t}{2} \sum_{i=1}^k y_i^2 + \left(\frac{(1+t\, n)}{2}  - \frac{(1-t)}{2}\right)\sum_{j=k+1}^ny_j^2\\
 \label{eq111}& \leq \dfrac{(n+1)\, t}{2}\|y\|^2.
\end{align}
Since $[(y,(u+\I_{rC}) (y)), (\Bar{y}, \bar{q}_t(\Bar{y}))]$ is tangent to $\bar{q}_t$ at $\Bar{y}$, we can apply \eqref{eq10.57} to $\Bar{y}$ instead of $x_0$, and obtain
\begin{align}\label{eq222}
\dist((y,(u+\I_{rC}) (y)), \epi (\bar{q}_t+ \I_C)) \geq \dfrac{t}{8}\|y-\bar{y}\|^2.    
\end{align}
Combining \eqref{eq111} and \eqref{eq222} with $y\in rC$, we deduce
$$\|y-\Bar{y}\|^2\leq \frac{8}{t}\dist((y,(u+\I_{rC})(y)), \epi (\bar{q}_t+ \I_C)) \leq 4(n+1)nr^2.$$

By the previous inequality, we deduce that
\begin{align*}
 \Bar{y}\in y+4\sqrt{n+1}rC \subset rC+4\sqrt{n+1}rC
  = (1+4\sqrt{n+1})\,rC \subset 4nrC.
\end{align*}
This completes the proof.
\end{proof}

\begin{lema} \label{lemaclaim3}
For $r>0$ sufficiently small,
\begin{align*}
 Z(u+\I_{rC})& \leq   \dfrac{\rho }{2}V_n(rC).
\end{align*}
\end{lema}

\begin{proof} 
By Lemma~\ref{claim3}, for sufficiently small $r>0$, the set $\conv(u+\I_{rC},\bar q_t+\I_C)$ differs from $\bar q_t+\I_C$ only in the neighborhood $4nrC$ of $x_0$.

Let $m_r$ be the maximal number of points $y_1, \dots, y_{m_r} \in C$ such that the sets 
$y_i + 8n r\, C$, $i = 1, \dots, m_r$, form part of a partition of $C$.
As before, we define
$$g_i(x+y_i)= (u+\I_{rC}+\ell_{y_i})(x),$$
and 
$$w_r=\conv(g_1, \dots,g_{m_r}, \bar{q}_t +\I_C),$$
where $\ell_{y_i}$ are the affine functions given by \eqref{psi} which satisfy \eqref{affine}. This construction implies that
$$w_r\seq \bar{q}_t+\I_C \quad \mbox{as} \quad r\rightarrow 0.$$
Moreover,
$$m_r\leq \dfrac{V_n(C)}{V_n( 8nr\, C)}$$
and
$$w_r(x+y_i)\leq (u+\I_{rC}+\ell_{y_i})(x)$$
for $i=1, \dots, m_r, x\in\mathbb{R}^n$, and $r>0$ sufficiently small. By construction,  the intersection of $y_i+rC$ and $y_j+rC$ for $i\neq j$ is either empty or at most an $(n-1)$-dimensional polytope. Since $Z$ is a  non-negative, simple, and  dually epi-translation invariant valuation, it follows that for every   $r>0$ sufficiently small, 
\begin{align*}
    Z(w_r) & \geq  m_r Z(u+\I_{rC}),
\end{align*}
and using  that $Z$ is  $\tau$-upper semicontinuous, we then obtain
\begin{align*}
    Z(\bar{q}_t+\I_C) & \geq \limsup_{r\rightarrow 0}Z(w_r) \\
    & \geq \limsup_{r\rightarrow 0}\dfrac{V_n(C)}{V_n( 8nr\, C)}  Z(u+\I_{rC}). 
\end{align*}
Thus, for every $\eta>0$,
$$\dfrac{V_n(C)}{V_n(8nr\, C)} Z(u+\I_{rC}) \leq Z(\bar{q}_t+\I_C) + \eta,$$
that is,
$$Z(u+\I_{rC}) \leq \dfrac{V_n(8nr\, C)}{V_n(C)}(Z(\bar{q}_t+\I_C)+ \eta) $$
for $r>0$ sufficiently small.

By \eqref{eq10.52}, using that $Z$ is dually-epi translation invariant and   applying the Claim, we find
\begin{align*}
    Z(u+\I_{rC}) &\leq \dfrac{V_n(8nr\, C)}{V_n(C)} \left( \dfrac{\rho \, V_n(C)}{4(4n)^n}+\eta \right)\\
    &\leq \dfrac{\rho}{4}V_n(rC)+\dfrac{\eta \, V_n(8nrC)}{V_n(C)}\\
    & =  V_n(rC)\left(\dfrac{\rho}{4}+\dfrac{\eta \, (8n)^n}{V_n(C)}\right).
\end{align*}
Choosing  $\eta = \dfrac{\rho \, V_n(C)}{4(8n)^n}$, we  obtain
\begin{align}\label{eq10.62}
 Z(u+\I_{rC})& \leq   \dfrac{\rho }{2}V_n(rC)
\end{align}
for $r>0$ sufficiently small.
\end{proof}

We now complete the proof of Proposition~\ref{regularr2}.

\begin{proof}[Proof of Proposition~\ref{regularr2}]
Let $r>0$ be sufficiently small so that the conclusion of Lemma~\ref{lemaclaim3} holds, and set
\[
P_r(x_0)=rC.
\]
Define
\[
v_r^{x_0}=l_{x_0}\vee l^c,
\]
where $l_{x_0}$ is the affine function parametrizing the tangent hyperplane to $u$ at $x_0$ and coinciding with $u$ at this point, while $l^c$ is the function chosen in \eqref{li_lc}. By construction,
\[
l^c(x)\le v_r^{x_0}(x)\le l^i(x)
\quad\text{for all }x\in\mathbb{R}^n,
\]
which gives \eqref{eq10.66m}.

On the other hand, by Lemma~\ref{lemaclaim3}, we have
\[
Z(u+\I_{P_r(x_0)})
\leq \frac{\rho}{2}V_n(P_r(x_0)),
\]
which gives \eqref{eq10.67m}. Finally, by \eqref{eq10.53.}, since $P_r(x_0)=rC$, we obtain
\[
\diam(P_r(x_0))^n
=
n^{\frac n2}V_n(P_r(x_0)).
\]
Thus there exists $r(x_0)>0$ such that \eqref{eq10.67m}, \eqref{eq10.66m}, \eqref{eq41112m}, and \eqref{eq10.53..m} hold for every $0<r\leq r(x_0)$. This completes the proof.
\end{proof}

\subsection{Proof of Proposition \ref{singular}}

Finally, we recall the statement of Proposition~\ref{singular}.

\begin{prop}\label{singularr}
Let $Z:\LP\RR\rightarrow\mathbb{R}$ be a simple, $\tau$-upper semicontinuous, equi-affine and  dually epi-translation   invariant valuation that   vanishes on  indicator functions of polytopes and cylinder functions. Let $u = \tilde{u}\,\square\, \frac{\lambda}{2}\big(\|\cdot\|^2+\I_{\mu C}\big)$, where $\tilde{u}\in \LP(\mathbb{R}^n)$,  $ \lambda>0$, and $0<\mu<1$.   There is a constant $c_{\lambda}$, depending on the dimension $n$ and on $\lambda$, such that
\begin{align*}
    Z(u+\I_{\dom u \, \cap \,P})\leq c_{\lambda}\, V_n(\dom u\cap P) 
\end{align*}
for every polytope $P\subset \mathbb{R}^n$.
\end{prop}

Let 
$$u = \tilde{u}\,\square\, \frac{\lambda}{2}\big(\|\cdot\|^{2}+\I_{\mu C}\big),
\qquad \tilde{u}\in \LP(\mathbb{R}^n),\ \lambda>0, \text{ and }\, 0<\mu< 1.$$
By Lemma \ref{gradient}, there exists $y_0\in \dom \tilde{u}$ such that the quadratic function
$$q_{x_0}(x) = \tilde{u}(y_0) + \frac{\lambda}{2}\|x - y_0\|^2$$
has the property that its graph touches the graph of $u$, restricted to $\dom u$, at the point $(x_0, u(x_0))$ with $x_0\in\dom u$. 
Moreover, $q_{x_0}(x)\ge u(x)$ for all $x\in y_0+\mu C$, and  also $x_0\in y_0+\mu C$. 

This implies that, for every sufficiently small $0<r<\mu$, we may choose a point
$p_r\in\interior(\dom u)$ such that
\begin{align}\label{pr}
x_0\in p_r+rC\subset y_0+\mu C,
\end{align}
and the epigraph of
\(
q_{x_0}+\I_{\dom u\cap(p_r+rC)}
\)
is contained in the epigraph of $u$. If  $x_0\in\interior(\dom u)$, we may simply take $p_r=x_0$. Otherwise, if $x_0\in \bd \dom u$, then $x_0\in \bd (p_r+rC)$. This second condition holds because $x_0\in y_0+\mu C$ and $y_0\in \dom \tilde{u}$, which for $x_0\in \bd \dom u$ implies that $x_0\in \bd (y_0+\mu C)$.

Let $\tau_{y_0}(x)=x-y_0$ for $x\in\mathbb{R}^n$. Note that 
$$q_\lambda(x)= q_{x_0}\circ\tau_{y_0}^{-1}(x) -\tilde{u}(y_0)= \dfrac{\lambda}{2}\sum_{i=1}^n x_i^2.$$
Moreover, $y_0+\mu C\subset y_0+ C$, which implies  $\tau_{y_0}(x_0)\in \tau_{y_0}(p_r+rC)\subset\tau_{y_0}(y_0+C)=C$. Without loss of generality, we may assume $\tilde{u}(y_0)=0$.

Arguing similarly to Lemma~\ref{claim1}, we show that the function
$$\conv\bigl(u\circ \tau_{y_0}^{-1}+\I_{\tau_{y_0}(\dom u\,\cap \, (p_r+r\, C))}, q_\lambda+\I_{C}\bigr)$$
coincides with $q_\lambda+\I_{C}$ outside a small neighborhood of $\tau_{y_0}(x_0)$.  Set 
$$v_r= u\circ \tau_{y_0}^{-1}+\I_{\tau_{y_0}(\dom u\,\cap \, (p_r+r\, C))}, \quad r>0.$$

\begin{lema}\label{claim4}
For every $x_0\in\dom u$, let $p_r=p_r(x_0)$ be chosen so as to satisfy \eqref{pr}. For $0<r<\mu$, define
$$
G_r = \left\{y\in\mathbb{R}^n\mid  \conv\bigl(v_r,\, q_\lambda+\I_{C}\bigr)(y) \neq (q_\lambda+\I_{C})(y) \right\}.
$$
Then 
$$
G_r \subset  p_r+\tau_{y_0}(3\sqrt{n} \, r\, C).
$$
\end{lema}

\begin{proof}
Let $y\in \tau_{y_0}(\dom u\,\cap \, (p_r+r\, C))$ and let  $(\bar{y}, q_{\lambda}(\bar{y}))$ be such that $[(\bar{y}, q_{\lambda}(\bar{y})), (y, v_r(y))]$ is tangent to $q_\lambda$ at $\bar{y}$. This implies 
$$\frac{2}{\lambda} v_r(y)+\|\bar{y}\|^2-2\langle \bar{y},y\rangle=0$$
and 
$$\|y-\bar{y}\|= \sqrt{\|y\|^2-\frac{2}{\lambda} v_r(y)}.$$ 
Moreover, since $(y, v_r(y))$ lies between the affine function associated with a subgradient of $v_r$ at $\tau_{y_0}(x_0)$ (which agrees with the gradient of $q_\lambda$ at $\tau_{y_0}(x_0)$) and the epi-graph of $q_\lambda$, we have
\begin{align*}
 \dist((y, v_r(y)), \epi(q_\lambda)) \leq \dist(\{y\}\times\mathbb{R}\cap H(\tau_{y_0}(x_0),v_r(\tau_{y_0}(x_0))), (y, q_\lambda(y))) 
 \leq   \frac{\lambda}{2} \|y-\tau_{y_0}(x_0)\|^2.
\end{align*}
Thus,
\begin{align*}
\|y-\bar{y}\|\leq \sqrt{\|y\|^2-\frac{2}{\lambda} v_r(y)} 
\leq \sqrt{\|y\|^2-\frac{2}{\lambda}\left(\frac{\lambda}{2}\|y\|^2-\frac{\lambda}{2}\|y-\tau_{y_0}(x_0)\|^2\right)}
=\|y-\tau_{y_0}(x_0)\|.
\end{align*}
We conclude that
\begin{align*}
G_r & \subset y+\, 2\sqrt{n}rC \subset p_r+\, \tau_{y_0}(rC)+2\sqrt{n}\, rC\subset p_r+ \tau_{y_0}(3\sqrt{n}\, rC).
\end{align*}
\end{proof}

\begin{lema}\label{localcube}
For every $x_0\in\dom u$, let $p_r=p_r(x_0)$ be chosen so as to satisfy \eqref{pr}. Then there exists
$r_0=r_0(\lambda,\mu)>0$ such that, for every $0<r<r_0$ and every sufficiently large $k$,
\begin{align*}
Z\bigl(u+\I_{\dom u\cap(p_r+(1+\frac1k)rC)}\bigr)
\leq
2^{3(n+1)}(\sqrt{n})^nZ(q_\lambda+\I_C)\, V_n(\dom u\cap (p_r+(1+\tfrac{1}{k})r\, C)).
\end{align*}
\end{lema}

\begin{proof}
Using Lemma~\ref{claim4}, we now construct convex functions $w_r$ that are
$\tau$-convergent to $q_\lambda+\I_C$. Let $m_r$ be the maximum number of points $y_1, \dots, y_{m_r}\in C$ such that the sets
\begin{align}\label{eq1.73}
y_i+(p_r+\tau_{y_0}(3\sqrt{n}\, rC))\cap C , \qquad i=1,\dots, m_r,    
\end{align}
form part of a partition of $C$. Consider the affine functions $\ell_{y_1}, \dots, \ell_{y_{m_r}}$, given by \eqref{psi}, associated to $y_i$ and $q_{\lambda}$, respectively. Define
$$g_i(y+y_i)=(v_r+\ell_{y_i})(y),$$
and
$$w_r= \conv\left(g_1+\I_{(-y_1+\tau_{y_0}(\dom u\, \cap \, (p_r+r\, C))}, \dots, g_{m_r}+\I_{(-y_{m_r}+\tau_{y_0}(\dom u\, \cap \, (p_r+r\, C))},q_\lambda+\I_{C}\right).$$
Then $w_r\seq q_\lambda+\I_{C}$ as $r\to 0$. Since $Z$ is a $\tau$-upper semicontinuous valuation, this implies that
\begin{align}\label{eq1.75}
  Z( q_\lambda+\I_{C})\geq \frac{1}{2}Z(w_r)  
\end{align}
for $0<r\leq r_1(\lambda)$ with some suitable $r_1(\lambda)>0$. 
\medskip

By \eqref{eq1.73} and  Lemma~\ref{claim4}, we have
$$w_r(y+y_i)\leq (v_r+\ell_{y_i})(y),$$
for $i=1, \dots, m_r$ and all $y\in\mathbb{R}^n$. 
Therefore, the intersection of $-y_i+\tau_{y_0}( \dom u\, \cap \, (p_r+r\, C))$ and $-
y_j+\tau_{y_0}(\dom u\, \cap \, (p_r+r\, C))$ for $i\neq j$ 
is either empty or an at most  $(n-1)$-dimensional polytope.
Since $Z$ is a non-negative, simple, and dually epi-translation invariant valuation, we get
\begin{align*}
 Z(w_r)& \geq \sum_{i=1}^{m_r}Z\bigl(v_r+\ell_{y_i}\bigr)= m_r\, Z\bigl(v_r\bigr). 
\end{align*}
It follows from this and  \eqref{eq1.75} that 
$$Z\bigl(v_r\bigr) \leq \frac{2}{m_r}Z(q_\lambda+ \I_{C}).$$
Since $m_r V_n((p_r+\tau_{y_0}(3\sqrt{n}r\, C))\cap C)\rightarrow V_n(C)$ as $r\rightarrow 0$,  there exists  $r_2>0$ such that
\begin{align}\label{eqr2}
 m_rV_n(3\sqrt{n}r\, C)\geq m_r V_n((p_r+\tau_{y_0}(3\sqrt{n}r\, C))\cap C)\geq \frac{1}{2}\, V_n(C)   
\end{align}
for $0<r\leq r_2$.
Hence,
\begin{align}\label{eq 76}
Z\bigl(v_r\bigr) \leq 4(3r)^n(\sqrt{n})^nZ(q_\lambda+\I_{C}).    
\end{align}

Recall that $\dom u= Q+\mu\, C$, where $Q=\dom \tilde{u}\in\mathcal{P}^n$ and $\mu>0$. Then 
\begin{align*}
   \dom u\cap (p_r+r\, C) &= (Q+\mu C)\cap (p_r+r\, C).
\end{align*}
Moreover, there are $z_1\in Q$ and $z_2\in C$ such that $p_r=z_1+\mu z_2$, which implies
\begin{align*}
    (z_1+\mu C)\cap (z_1+\mu z_2+r\, C)\subset (Q+\mu C)\cap (z_1+\mu z_2+r\, C), 
\end{align*}
and consequently,
\begin{align}
\nonumber   V_n(\dom u\cap (p_r+r\, C))& \geq V_n((z_1+\mu C)\cap (z_1+\mu z_2+r\, C))\\
\nonumber & = V_n(\mu C\cap (\mu z_2+r\, C))\\
\label{eq 77}    &\geq \frac{r^n}{2^n} 
\end{align}
for $0<r\leq r_3(\mu)$. Using this in \eqref{eq 76}, and the vertical and horizontal translation invariance of $Z$, and \eqref{mudanca_zeta}, we obtain
\begin{align}\label{cube}
Z(u+\I_{ \dom u\,\cap \,(p_r+rC)}) =Z\bigl(v_r\bigr) \leq  2^{3n+2}(\sqrt{n})^nZ(q_\lambda+\I_C)\, V_n(\dom u\cap (p_r+r\, C))
\end{align}
for $0<r< \min\left\{r_1(\lambda), r_2, r_3(\mu)\right\}$.

Let $x_0\in  \dom u$ and let $0<r< \min\left\{r_1(\lambda),r_2, r_3(\mu)\right\}$ be such that \eqref{cube} holds. For each $k\in\mathbb{N}$ define
$$v_{r,k}=  u\circ \tau_{y_0}^{-1}+\I_{\tau_{y_0}(\dom u\, \cap \, (p_r+(1+\frac{1}{k})r\, C))}.$$
By construction $v_{r,k}$ is $\tau$-convergent to $v_r$ as $k\to \infty$. Since $Z$ is $\tau$-upper semicontinuous we obtain 
$$
\limsup_{k\to\infty} Z(v_{r,k}) \le Z(v_r).
$$
Hence, for any fixed $\varepsilon>0$ there exists $k_0$ such that for all $k\ge k_0$,
$$
Z(v_{r,k}) \le Z(v_r) + \varepsilon.
$$
Combining this with the estimate previously obtained for $Z(v_r)$ yields that for $k$ large enough
\begin{align*}
\nonumber Z(u+\I_{\dom u\, \cap \, (p_r+(1+\frac{1}{k})r\, C)})&=Z(v_{r,k})\\
\nonumber &\le Z(v_r)+  2^{3n+2}(\sqrt{n})^nZ(q_\lambda+\I_C)\, V_n(\dom u\cap (p_r+r\, C))\\
&\le  2^{3(n+1)}(\sqrt{n})^nZ(q_\lambda+\I_C)\, V_n(\dom u\cap (p_r+(1+\tfrac{1}{k})r\, C)).
\end{align*}
This completes the proof.
\end{proof}

Note that $x_0\in \dom u$ is an interior point of $p_r+(1+\frac{1}{k})r\, C$, since the cube has been strictly enlarged from $r$ to $(1+\frac{1}{k})r$.

We now complete the proof of Proposition~\ref{singularr} using Lemma~\ref{localcube}.

\begin{proof}[Proof of Proposition \ref{singularr}]
Let $P\subset\mathbb{R}^n$ be an arbitrary polytope. Choose an open set $U$ such that
\begin{align}\label{eq 78}
\dom u\cap P  \subset U,    
\end{align}
and
\begin{align}\label{eq 79}
V_n(U)\leq 2V_n(\dom u\cap P).    
\end{align}

Let $\mathcal{C}$ be the family of closed cubes $\bar{C}=\bar{C}(x,\bar{r})$ with center $x\in \interior(\dom u) \cap P$ and side length $2\bar{r}$, where
$\bar{r}=\Big(1+\tfrac{1}{k_x}\Big)r$ and $0<r< \frac{1}{(1+2\sqrt{n})}r_0$ for some $k_x\in\mathbb{N}$  chosen so that the estimate in Lemma~\ref{localcube} applies,  and such that $\dom u\cap  \bar{C} \subset U$. Then the  interiors of the cubes  $\bar{C}$ with $\bar{C}\in \mathcal{C}$ form an open cover of $\dom u \cap P$, since by the construction we know that  every boundary point of 
$\dom u$ is contained as an interior point of some cube of the form $x_r+\bigl(1+\tfrac{1}{k}\bigr)rC$. Since $\dom u\, \cap \,P$ is compact, we can select a finite subcover,  denoted by $\mathcal{L}\subset \mathcal{C}$.  
By  a standard argument known as Vitali's Lemma (cf.\ \cite[Lemma 1.9]{falconer1985geometry}),  we may choose pairwise disjoint cubes $\bar{C}(x_1, \bar{r}_1), \dots, \bar{C}(x_m, \bar{r}_m)$ from $\mathcal{C}$  such that
$$\dom u \cap P \subset D=\bigcup_{i=1}^m \bar{C}\left(x_i, (1+2\sqrt{n})\, \bar{r}_i\right).$$
Since $Z$ is non-negative,
\begin{align}
 Z\bigl(u+\I_{  \dom u\, \cap \, D}\bigr)& \le \sum_{i=1}^m Z\bigl(u+ \I_{ \dom u\, \cap \, \bar{C}(x_i, (1+2\sqrt{n})\,\bar{r}_i)}\bigr),
\end{align}
and by Lemma \ref{localcube}, it follows that 
\begin{align}
\nonumber Z\bigl(u+\I_{  \dom u\, \cap \, D}\bigr) & \leq  2^{3(n+1)}(\sqrt{n})^n(1+2\sqrt{n})^n\,Z(q_\lambda+ \I_C)\sum_{i=1}^mV_n(\bar{r}_i\, C)\\
\nonumber & =   2^{5n+4}(\sqrt{n})^n(1+2\sqrt{n})^n\,Z(q_\lambda+ \I_C) \sum_{i=1}^m \frac{1}{2}\Big(\frac{{\bar{r}_i}}{2}\Big)^n.
\end{align}
Since the cubes $\bar{C}(x_i,\bar{r}_i)$ are pairwise disjoint, we obtain by \eqref{eq 77} that
$$\sum_{i=1}^m \frac{1}{2}\Big(\frac{{\bar{r}_i}}{2}\Big)^n\leq \frac{1}{2}\sum_{i=1}^m V_n(\dom u\cap \bar{C}(x_i,\bar{r}_i)),$$
and by \eqref{eq 78} and  \eqref{eq 79} that
$$\frac{1}{2}\sum_{i=1}^m V_n(\dom u\cap \bar{C}(x_i,\bar{r}_i))\leq \frac{1}{2}V_n(U)\le 
V_n(\dom u\cap P).$$

Finally, since $\dom u\cap P\subset D$ and $Z$ is non-negative and simple, we obtain
\[
Z(u+\I_{\dom u\cap P})
\leq
Z(u+\I_{\dom u\cap D})
\leq
c_{\lambda}\,V_n(\dom u\cap P),
\]
where
\(
c_{\lambda}
=
2^{5n+4}(\sqrt n)^n(1+2\sqrt n)^n Z(q_\lambda+\I_C)\).
This completes the proof.
\end{proof}

\end{document}